\documentstyle[12pt]{article}
\input amssym.def
\begin {document}
\def \Z{{\bf Z}}
\def \C{{\bf C}}

\def \N{\Bbb N}

\def \0{{\bf 0}}
\def \x{{\bf x}}
\def \y{{\bf y}}
\def \z{{\bf z}}
\def \k{{\bf k}}
\def \m{{\bf m}}

\def \d{{\bf d}}
\def \D{{\bf D}}
\def \r{{\bf r}}
\def \s{{\bf s}}

\def \fg{\frak g}

\def \Res{{\rm Res}}
\def \End{{\rm End}\;}

\def \Hom{{\rm Hom}}

\def \<{\langle} 
\def \>{\rangle} 

\def \a{\alpha }

\def \be{\begin{equation}\label}
\def \ee{\end{equation}}
\def \bexa{\begin{exa}\label}
\def \eexa{\end{exa}}

\def \bl{\begin{lem}\label}
\def \el{\end{lem}}
\def \bt{\begin{thm}\label}
\def \et{\end{thm}}
\def \bp{\begin{prop}\label}
\def \ep{\end{prop}}
\def \br{\begin{rem}\label}
\def \er{\end{rem}}
\def \bc{\begin{coro}\label}
\def \ec{\end{coro}}
\def \bd{\begin{de}\label}
\def \ed{\end{de}}
\def \pf{{\bf Proof. }}

\newtheorem{thm}{Theorem}[section]
\newtheorem{prop}[thm]{Proposition}
\newtheorem{coro}[thm]{Corollary}

\newtheorem{exa}[thm]{Example}
\newtheorem{lem}[thm]{Lemma}
\newtheorem{rem}[thm]{Remark}
\newtheorem{de}[thm]{Definition}

\makeatletter
\@addtoreset{equation}{section}

\makeatother
\makeatletter

\baselineskip=16pt
\begin{center}{\Large \bf On certain
higher dimensional analogues of vertex algebras}

\vspace{0.5cm}
Haisheng Li\footnote{Partially supported by NSF grant
DMS-9970496 and a Rutgers Research Council grant}\\ 
Department of Mathematical Sciences, Rutgers University,
Camden, NJ 08102\\
and\\
Department of Mathematics, Harbin Normal University,
Harbin, P.R. China
\end{center}

\begin{abstract} 
A higher dimensional analogue of the notion of vertex  algebra 
is formulated in terms of formal variable language with 
Borcherds' notion of $G$-vertex algebra as a motivation.
Some examples are given and certain analogous duality properties 
are proved. Furthermore, it is proved that for any vector space $W$, 
any set of mutually local
multi-variable vertex operators on $W$ in a certain canonical way 
generates a vertex algebra with $W$ as a natural module. 
\end{abstract}

\section{Introduction}

Vertex (operator) algebras, introduced in mathematics ([B1], [FLM]),
are known essentially to be chiral algebras of
two dimensional conformal field theories, as formulated
in physics ([BPZ], [MS]). 
In [B2], higher dimensional analogues of vertex algebras were 
established in terms of what were called relaxed multi-linear categories.
(In fact, the setting of [B2] is extremely general.)
It was proved in [Sn] (cf. [B2]) that this notion in a special case 
is equivalent to the notion of vertex algebra defined in [B1] and [FLM] 
in terms of formal variables. It is hoped that these algebras have 
the same relation to higher dimensional
quantum field theories that vertex algebras have to one dimensional 
quantum field theories (or to ``chiral halves'' of two dimensional 
conformal field theories). 

In vertex algebra theory, even though the notion of ordinary vertex 
algebra can be formulated in terms of ${\cal{D}}$-modules
([BD], [HL]) and in terms of relaxed multi-linear categories ([B2], [Sn]),
the formal variable approach plays a unique role.
Analogously, a formal variable formulation of
higher dimensional analogues of vertex algebras should  be of importance.
The main purpose of this paper is to formulate certain higher dimensional
analogues of vertex algebras in terms of formal variable language in 
the same way that the ordinary vertex algebras were defined in [B1] and [FLM].

Vertex algebras are often thought of as a kind of generalized algebras 
equipped with infinitely many bilinear operations.
Another point of view is that
vertex operator algebras are ``algebras'' of vertex operators just as
classical (associative, or Lie) algebras are algebras of linear
operators, where vertex operators on a vector space $W$ are understood 
as elements of the space $\Hom (W,W((x)))$.
(This particular viewpoint was emphasized in [FLM].)
Then for a vertex algebra $V$, each element $v$ is represented
by the vertex operator $Y(v,x)$ on $V$, where $Y(v,x)$ is
an analogue of the left multiplication associated to an element
for a classical non-associative algebra.
In the definition of the notion of vertex algebra, 
the main axiom is what was called in [FLM]
and [FHL] the Jacobi identity:
\begin{eqnarray}
& &z^{-1}\delta\left(\frac{x-y}{z}\right)Y(u,x)Y(v,y)w
-z^{-1}\delta\left(\frac{-y+x}{z}\right)Y(v,y)Y(u,x)w\nonumber\\
&=&y^{-1}\delta\left(\frac{x-z}{y}\right)Y(Y(u,z)v,y)w.
\end{eqnarray}

For higher dimensions, we again use the viewpoint that
vertex (operator) algebras are ``algebras'' of local {\em multi-variable}
vertex operators. Given a vector space $W$, we naturally consider
elements of the space
\begin{eqnarray}
\Hom \left(W,W[[x_{1},\dots,x_{n}]][x_{1}^{-1},\dots,x_{n}^{-1}]\right)
\end{eqnarray}
as our $n$-variable vertex operators. For an $n$-dimensional analogue
$V$ of vertex algebras, each element $v$ of $V$ should be represented
by an $n$-variable vertex operator $Y(v,\x)$ on $V$, where 
$\x=(x_{1},\dots,x_{n})$. 
As in the ordinary case [DL] (cf. [FHL], [FLM], [G]), 
vertex operators $Y(u,\x)$ and $Y(v,\y)$
are assumed to be mutually local in the sense that 
\begin{eqnarray}
(\x-\y)^{k}Y(u,\x)Y(v,\y)=(\x-\y)^{k}Y(v,\y)Y(u,\x)
\end{eqnarray}
for some $k\in \N$. (See Section 2 for notations.) Note that
in the 1-dimensional case , product $Y(u,x)Y(v,y)$ 
is described by rational functions with only three possible 
poles at $x=0,\infty, y$ and the Jacobi identity is a version of 
classical Cauchy residue theorem ([FLM], [FHL]).
In the $n$-dimensional case, from locality product $Y(u,\x)Y(v,\y)$ 
can be described by rational functions with only $3^{n}$ possible poles.
But, there are 
still three {\em explicit} relevant quantities (two products and one iterate).
What we do next is to employ formal Laurent series expansions 
in certain domains.
For each variable $x_{i}$, we have two disjoint domains
$|x_{i}|>|y_{i}|>0$ and $|y_{i}|>|x_{i}|>0$. This gives us
total $2^{n}$ different domains for $\x=(x_{1},\dots,x_{n})$.
Each such domain can be represented by a binary code 
$\alpha\in (F_{2})^{n}$ of length $n$ in the obvious way.
For $\m\in \Z^{n}$, we then define $2^{n}$ formal series 
$\iota_{\alpha,\x,\y}((\x-\y)^{\m}Y(u,\x)Y(v,\y))$
which roughly speaking are the analytic continuation of 
$(\x-\y)^{\m}Y(u,\x)Y(v,\y)$ in the $2^{n}$ different domains.
Note that the products $(\x-\y)^{\m}Y(u,\x)Y(v,\y)$ and 
$(-y+\x)^{\m}Y(v,\y)Y(u,\x)$ are two special cases.
So, each $\iota_{\alpha,\x,\y}((\x-\y)^{\m}Y(u,\x)Y(v,\y))$ 
can be considered as an implicit product.
Now our multi-variable Jacobi identity reads as
\begin{eqnarray}
& &\sum_{\alpha\in (F_{2})^{n}}(-1)^{|\alpha|}\iota_{\alpha,\x,\y}
\left(\z^{-1}\delta\left(\frac{\x-\y}{\z}\right)
Y(u,\x)Y(v,\y)\right)\nonumber\\
&=&\y^{-1}\delta\left(\frac{\x-\z}{\y}\right)Y(Y(u,\z)v,\y),
\end{eqnarray}
where $\x,\y,\z$ are $n$-dimensional variables and 
the delta functions are $n$-dimensional.
In the $1$-dimensional case, Jacobi identity 
does not require locality, instead, it implies locality.
For higher dimensional cases, to state the Jacobi identity,
or for $\iota_{\alpha,\x,\y}((\x-\y)^{\m}Y(u,\x)Y(v,\y))$ 
to be defined, we need the locality of $Y(u,\x)$ and $Y(v,\y)$
as a prerequisite. 
Now, the main axioms of 
our higher dimensional analogue of the notion of vertex algebra
are locality and the above Jacobi identity.
Typical examples of such vertex algebras are tensor products of
$n$ ordinary vertex algebras.

Once this higher dimensional analogue of the notion of 
vertex algebra is defined, all the results in the ordinary case 
(see [FLM], [FHL], [DL], [Li2]) are carried over, and 
what is more, the same old proofs work fine almost all the time.
Especially, it is proved that locality (or weak commutativity) 
together with some minor information is strong enough to recover 
the Jacobi identity.
Thus, one may define the notion of vertex algebra by using
locality instead of the Jacobi identity as the main axiom
(see [DL], [FHL], [FLM], [G], [Li2], see also [Ka]).
This is analogous to a phenomenon in the classical commutative associative 
algebra theory. It is known that
vertex algebras are analogous to commutativity 
associative algebras with identity while locality is analogous  
to commutativity
for left multiplications associated with elements of 
classical algebras. As it was explained in [Li2],
commutativity associative algebras
with identity are exactly non-associative algebras with identity
with all left multiplications being commuting.
However, in defining the notion of module for a commutative
associative algebra we use {\em associativity} instead of
{\em commutativity}. It seems that it is a universal principle 
that the notion of module should be defined by using all 
the axioms that make sense in defining the notion of the algebra.
In view of this, using Jacobi identity to define
the notion of vertex algebra is more natural than using locality.

As a good practice of the viewpoint that vertex (operator) algebras 
are ``algebras'' of vertex operators, in [Li2] we established
an analogue of endomorphism algebra.
It was proved that any maximal space of pair-wise mutually local 
vertex operators on a vector space $W$ has a canonical vertex algebra 
structure with $W$ as a natural module. In terms of vertex algebras of 
local vertex operators, the notion of module and the notion of 
representation are identified just as in classical theories.
Furthermore, it follows from Zorn's Lemma that 
for an abstract vector space $W$,
any set of mutually local vertex operators on $W$ 
in a certain canonical way generates a vertex algebra
(a subspace of $\Hom (W,W((x)))$) with $W$ as a 
natural module. As mentioned in [Li2], this result is 
an analogue of the trivial classical fact --- any set of commuting 
linear operators on a vector space $U$
generates a commutative associative algebra with $U$ as a module.
In addition to its conceptual importance, this result is 
quite useful to construct examples of vertex algebras 
and their modules. 

Certain results similar to that of [Li2] were also independently 
obtained in [LZ] and the work [LZ] is also in this spirit. 

In this paper, we establish the main results of [Li2] for higher 
dimensional analogues of  vertex algebras. In fact, this paper is 
modeled on [Li2]. In the 1-dimensional case, the adjoint vertex operator 
(vertex operators on the space of vertex operators) was defined 
in terms of commutators (or more generally cross brackets) and 
normal-order products. 
For higher dimensions, normal-order product is not defined 
in general. Note that 
for mutually local vertex operators $\psi,\phi$ and for $\m\in \Z^{n}$,
even though
$$ (\psi_{[\m]}\phi)(\y)=
\Res_{\x}((\x-\y)^{\m}\psi(\x)\phi(\y)-(-\y+\x)^{\m}\phi(\y)\psi(\x))$$
defines  a vertex operator, it is not what we want.
(In general, $\psi_{[\m]}\phi$ may not be local with the old vertex
operators $\psi$ and $\phi$.)
Here, we use an analogue of the Jacobi identity 
(or just the iterate formula) to
define the adjoint vertex operators.
Then the main results of [Li2] for higher dimensions are carried over.



This paper is motivated by [B2] and presumably, 
the $n$ dimensional analogue formulated in this paper is equivalent 
to Borcherds' notion in [B2] of commutative $G_{n}$-vertex algebra,
where $G_{n}=(H_{n},K_{n})$ is what was called therein
an elementary vertex group with
\begin{eqnarray}
H_{n}&=&{\C}[D_{1},\dots, D_{n}],\\
K_{n}&=&{\C}[[x_{1},\dots,x_{n}]][x_{1}^{-1},\dots,x_{n}^{-1}].
\end{eqnarray}
(See Remark 2.5 for more information.) 
Most likely, Snydal's proof in [Sn] will extend to higher dimensions
with some changes. In view of this, our $n$ dimensional analogues of 
vertex algebras are named as $G_{n}$-vertex algebras and the vertex 
operators we consider are named as $G_{n}$-vertex operators..

In [KO], a certain notion of vertex algebra (with $n=2$) was introduced
with an axiom in terms of operator product expansions. As mentioned in [KO],
if we only allow integral powers of formal variables in 
the definition of the notion of vertex algebra in [KO],
locality follows from the axioms. Then as mentioned before,
the Jacobi identity will follow. Therefore,
the notion of vertex algebra defined in [KO] in the special case
is equivalent to the notion of $G_{2}$-vertex algebra defined in this paper.
On the other hand, we emphasize that
our definition of vertex algebra naturally give rise to 
the (right) definition of module and that the main axiom, Jacobi identity,
is crucial in constructing vertex algebras from given vertex operators.

In the ordinary case, among the important and interesting
vertex (operator) algebras are those constructed from certain 
highest weight representations of affine Lie algebras (see [FZ]).
Then one naturally wants to consider the multi-variable 
generalizations of affine Lie 
algebras --- toroidal Lie algebras (cf. [BBS]). 
For toroidal Lie algebras, one can form generating functions 
in $n$-variables and from [BBS], one can easily check that all 
generating series are pair-wise mutually local. 
However, on the modules $W$ constructed in [BBS], generating series
$a(\x)$ for $a\in \fg$ do not give rise to $G_{n}$-vertex operators on $W$ 
because $a(\x)v$ in general involves infinitely many negative powers of
variables $x_{2},\dots,x_{n}$. So, our results here cannot apply to
toroidal Lie algebras. 
We believe that it is possible to relate this type of modules for 
toroidal Lie algebras to Borcherds' $G$-vertex algebras of a certain type.
This is currently under investigation in [Li3]. 	

This paper is organized as follows: In Section 2, we define 
certain generalized iota-maps and present some basic properties.
In Section 3, we define higher-dimensional analogues of vertex algebras
and present basic duality properties. In Section 4,
we present the theory of vertex algebras of multi-variable vertex operators.

\section{Basics of formal calculus and generalized iota-maps}
In this section, we first collect basic notations and conventions 
in formal calculus and we then define $G_{n}$-vertex operators 
and certain generalized iota-maps.
We also present some basic properties of the generalized iota-maps.

We here mainly follow [FLM] and [FHL] for the treatment of formal calculus. 
Throughout this paper, 
$x, x_{0}, x_{1}, x_{2},\dots, y, y_{0}, y_{1},y_{2},\dots, 
z, z_{0}, z_{1},z_{2},\dots$
are independent commuting  formal  variables. We typically use $U$ for a 
general vector space in establishing certain definitions and notations.

Let $n$ be a positive integer fixed throughout this paper. 
Set
\begin{eqnarray}
\x=(x_{1},\dots,x_{n}),\;\;\y=(y_{1},\dots,y_{n}),\;\;
\z=(z_{1},\dots,z_{n}).
\end{eqnarray}
For $\m=(m_{1},\dots,m_{n})\in {\Z}^{n}$, we set
\begin{eqnarray}
& &\x^{\m}=x_{1}^{m_{1}}\cdots x_{n}^{m_{n}}.
\end{eqnarray}
For $k\in {\Z}$, set
\begin{eqnarray}
\x^{k}=x_{1}^{k}\cdots x_{n}^{k}.
\end{eqnarray}
Then, for $\m\in \Z^{n},\;k\in \Z$, we write
$\x^{\m+k}$ for $\x^{\m}\x^{k}$. 
In particular, $\x^{-\m-1}$ is defined this way. 

For a vector space $U$,
denote by $U[[x_{1},x_{1}^{-1},\dots, x_{n},x_{n}^{-1}]]$ 
the space of all formal integral power series in $x_{1},\dots,x_{n}$ 
with coefficients in $U$.
That is,
\begin{eqnarray}
U[[x_{1},x_{1}^{-1},\dots, x_{n},x_{n}^{-1}]]
=\left\{ \sum_{\m\in \Z^{n}}f(\m)\x^{\m}\;|\; f(\m)\in U
\;\;\;\mbox{ for }\m\in \Z^{n}\right\}.
\end{eqnarray}
Denote by $U[[x_{1},\dots,x_{n}]]$
the subspace of all formal nonnegative power series  and
by $U((x_{1},\dots, x_{n}))$ the subspace of all formal 
Laurent (lower truncated) series, i.e.,
\begin{eqnarray}
U((x_{1},\dots,x_{n}))=U[[x_{1},\dots,x_{n}]][x_{1}^{-1},\dots,x_{n}^{-1}].
\end{eqnarray}
We shall also use the notations $U[[\x]]$ and $U((\x))$ for
$U[[x_{1},\dots,x_{n}]]$ and $U((x_{1},\dots,x_{n}))$, respectively.

In formal calculus, we shall heavily use the following binomial 
expansion convention:
\begin{eqnarray}
(x-y)^{m}=\sum_{i\ge 0}{m\choose i}(-1)^{i}x^{m-i}y^{i}\in \C[x,x^{-1}][[y]]
\end{eqnarray}
for $m\in \Z$, where ${m\choose i}={1\over i!}m(m-1)\cdots (m-i+1)$.

For $k\in {\Z}$, $\m=(m_{1},\dots,m_{n})\in {\Z}^{n}$, we define
\begin{eqnarray}
& &(\x-\y)^{k}=(x_{1}-y_{1})^{k}\cdots (x_{n}-y_{n})^{k},\\
& &(\x-\y)^{\m}=(x_{1}-y_{1})^{m_{1}}\cdots (x_{n}-y_{n})^{m_{n}}.
\end{eqnarray}
Then, for $\m\in \Z^{n},\;k\in \Z$, we define
$(\x-\y)^{\m+k}$ to be $(\x-\y)^{\m}(\x-\y)^{k}$.
For $\m=(m_{1},\dots,m_{n})\in {\Z}^{n}$ and $\i=(i_{1},\dots,i_{n})\in \N^{n}$, 
we set
\begin{eqnarray}
{\m\choose \i}={m_{1}\choose i_{1}}\cdots {m_{n}\choose i_{n}}.
\end{eqnarray}
Then
\begin{eqnarray}	
(\x+\y)^{\m}=\sum_{\i\in \N^{n}}{\m\choose \i}\x^{\m-\i}\y^{\i}.
\end{eqnarray}

Recall from [FLM] that for $A(\x),
B(x), C(\x)\in (\End U)[[x_{1},x_{1}^{-1},\dots, x_{n},x_{n}^{-1}]]$, 
the following associativity
\begin{eqnarray}
A(\x)(B(\x)C(\x))=(A(\x)B(\x))C(\x)
\end{eqnarray}
holds {\em provided} that all the products $A(\x)B(\x)C(\x)$, 
$A(\x)B(\x)$ and $B(\x)C(\x)$ exist.
The assumption of the existence of the products is crucial.
As it was illustrated in [FLM], abuse of this associativity could 
lead to paradoxes. For convenience, we refer 
the {\em associativity assumption} for the product $A(\x)B(\x)C(\x)$ to
the existence of the three products.

\br{rcancellationfact}
{\em We have the following cancelation law. Let 
$$f(\x,\y),\; g(\x,\y)\in 
U[[x_{i},x_{i}^{-1},y_{i},y_{i}^{-1}\;|\; i=1,\dots,n]]$$
be such that $f(\x,\y)$ and $g(\x,\y)$ involve only finitely many
negative powers of $x_{i}$ for $i\in I$,
where $I$ is a subset of $\{ 1,\dots,n\}$.
Suppose that there exist nonnegative integers $k_{i}$ for $i\in I$
such that
\begin{eqnarray}\label{epf=pg}
\left(\coprod_{i\in I}(x_{i}-y_{i})^{k_{i}}\right)f(\x,\y)
=\left(\coprod_{i\in I}(x_{i}-y_{i})^{k_{i}}\right)g(\x,\y).
\end{eqnarray}
Then $f(\x,\y)=g(\x,\y)$.
Indeed, under the given assumption on $f(\x,\y)$ and $g(\x,\y)$,
the associativity assumption holds for the products
$$\left(\coprod_{i\in I}(-y_{i}+x_{i})^{-k_{i}}\right)
\left(\coprod_{i\in I}(x_{i}-y_{i})^{k_{i}}\right)f(\x,\y)$$
and
$$\left(\coprod_{i\in I}(-y_{i}+x_{i})^{-k_{i}}\right)
\left(\coprod_{i\in I}(x_{i}-y_{i})^{k_{i}}\right)g(\x,\y).$$
Then multiplying both sides of (\ref{epf=pg}) by
$\coprod_{i\in I}(-y_{i}+x_{i})^{-k_{i}}$ and using associativity
we obtain $f(\x,\y)=g(\x,\y)$.}
\er

The formal residue $\Res_{x}f(x)$ of a formal series 
$f(x)=\sum_{n\in {\Z}}f(n)x^{n}\in U[[x,x^{-1}]]$ is defined to be 
the coefficient $f(-1)$ of $x^{-1}$ in $f(x)$.
We shall use the following convention:
\begin{eqnarray}
\Res_{\x}=\Res_{x_{1}}\cdots \Res_{x_{n}}.
\end{eqnarray}

The formal delta function is defined to be a formal series:
\begin{eqnarray}
\delta(x)=\sum_{m\in \Z}x^{m}\in {\C}[[x,x^{-1}]].
\end{eqnarray}
It is a so-called expansion of zero (see [FLM]):
\begin{eqnarray}
\delta(x)=(1-x)^{-1}-(-x+1)^{-1}
\end{eqnarray}
in terms of the binomial expansion convention.
Basic properties of delta function are:
\begin{eqnarray}
& &f(z)\delta(z)=f(1)\delta(z)\;\;\;\mbox{ for }f(z)\in {\C}[z,z^{-1}],\\
& &g(z_{1},z_{2})z_{2}^{-1}\delta\left({z_{1}\over z_{2}}\right)=
g(z_{1},z_{1})z_{2}^{-1}\delta\left({z_{1}\over z_{2}}\right)
\end{eqnarray}
for $g(z_{1},z_{2})\in {\C}[[z_{1},z_{1}^{-1},z_{2},z_{2}^{-1}]]$ 
such that $g(z_{1},z_{1})$ exists (algebraically), e.g., 
$g(z_{1},z_{2})\in {\C}((z_{1},z_{2}))$. Furthermore,
\begin{eqnarray}\label{edeltathreesub}
h(z_{0},z_{1},z_{2})z_{0}^{-1}\delta\left(\frac{z_{1}-z_{2}}{z_{0}}\right)
=h(z_{1}-z_{2},z_{1},z_{2})
z_{0}^{-1}\delta\left(\frac{z_{1}-z_{2}}{z_{0}}\right)
\end{eqnarray}
for $h(z_{0},z_{1},z_{2})\in U((z_{0},z_{1},z_{2}))$.

The following are the fundamental delta function identities 
in the formal variable approach to vertex (operator) algebras:
\begin{eqnarray}
& &z_{0}^{-1}\delta\left(\frac{z_{1}-z_{2}}{z_{0}}\right)
=z_{1}^{-1}\delta\left(\frac{z_{0}+z_{2}}{z_{1}}\right),\label{ethreechange}\\
& &z_{0}^{-1}\delta\left(\frac{z_{1}-z_{2}}{z_{0}}\right)
-z_{0}^{-1}\delta\left(\frac{-z_{2}+z_{1}}{z_{0}}\right)
=z_{2}^{-1}\delta\left(\frac{z_{1}-z_{0}}{z_{2}}\right).\label{ethreeidentity}
\end{eqnarray}

We define the following multi-variable delta functions:
\begin{eqnarray}
\delta(\x)&=&\delta(x_{1})\cdots \delta(x_{n}),\\
\y^{-1}\delta\left({\x\over \y}\right)
&=&y_{1}^{-1}\delta\left({x_{1}\over y_{1}}\right)
\cdots y_{n}^{-1}\delta\left({x_{n}\over y_{n}}\right),\\
\z^{-1}\delta\left(\frac{\x-\y}{\z}\right)&=&
z_{1}^{-1}\delta\left(\frac{x_{1}-y_{1}}{z_{1}}\right)\cdots 
z_{n}^{-1}\delta\left(\frac{x_{n}-y_{n}}{z_{n}}\right).
\end{eqnarray}
Then using (\ref{ethreechange}) we have
\begin{eqnarray}\label{exzchange}
\z^{-1}\delta\left(\frac{\x-\y}{\z}\right)
=\x^{-1}\delta\left(\frac{\z+\y}{\x}\right).
\end{eqnarray}
Furthermore, using (\ref{edeltathreesub}) we get
\begin{eqnarray}\label{esubstitution}
\z^{-1}\delta\left(\frac{\x-\y}{\z}\right)f(\x,\y,\z)
=\z^{-1}\delta\left(\frac{\x-\y}{\z}\right)f(\x,\y, \x-\y)
\end{eqnarray}
for $f(\x,\y,\z)\in U((\x,\y,\z))$.

Now we discuss a multi-dimensional analogue of (\ref{ethreeidentity}).
First,  (\ref{ethreeidentity}) immediately gives
\begin{eqnarray}\label{eformalproduct}
\prod_{i=1}^{n}\left(z_{i}^{-1}\delta\left(\frac{x_{i}-y_{i}}{z_{i}}\right)
-z_{i}^{-1}\delta\left(\frac{-y_{i}+x_{i}}{z_{i}}\right)\right)
=\prod_{i=1}^{n}y_{i}^{-1}\delta\left(\frac{x_{i}-z_{i}}{y_{i}}\right)
=\y^{-1}\delta\left(\frac{\x-\z}{\y}\right).
\end{eqnarray}
To express the left-hand side
in terms of $n$-dimensional delta functions, we need to expand
it into $2^{n}$-terms. 

Let $F_{2}=\{0,1\}$ be the 2-element field. Then 
$(F_{2})^{n}$ is an $n$-dimensional vector space over $F_{2}$.
For $\alpha\in (F_{2})^{n}$, denote by $|\alpha|$ 
the weight of $\alpha$, which is the number of 
nonzero $\alpha_{i}$'s in the expression 
$\alpha=(\alpha_{1},\dots,\alpha_{n})$. Denote by $e_{i}$ the element of 
$(F_{2})^{n}$
all of whose coordinates are $0$ except the $i$th coordinate which is $1$.
For $1\le i\le n$, we define a translation $\sigma_{i}$ 
on $(F_{2})^{n}$ by
\begin{eqnarray}
\sigma_{i}(\alpha)=\alpha+e_{i}\;\;\;\mbox{ for }\alpha\in (F_{2})^{n}.
\end{eqnarray}
Then $\sigma_{i}^{2}=1$ and
\begin{eqnarray}
(-1)^{|\sigma_{i}(\alpha)|}=-(-1)^{|\alpha|}.
\end{eqnarray}
Set
\begin{eqnarray}
E_{n}=\{\alpha\in (F_{2})^{n}\;|\; |\alpha|\in 2\Z\}.
\end{eqnarray}
Then
\begin{eqnarray}
(F_{2})^{n}=E_{n}\cup \sigma_{i}(E_{n})\;\;\;\mbox{(disjoint union)}
\end{eqnarray}
for $i=1,\dots, n$.

In the following we shall extend the iota-maps defined in [FLM] and [FHL]
and in [Le] to define linear maps $\iota_{\alpha,\x,\y}$ for 
$\alpha\in (F_{2})^{n}$. Their common domain will be 
the following subalgebra of the field of fractions of 
the ring $\C[[\x,\y]]$:
\begin{eqnarray}
A=\C[[\x,\y]][x_{i}^{-1},y_{i}^{-1}, (x_{i}+y_{i})^{-1}, (x_{i}-y_{i})^{-1}
\;|\; i=1,\dots,n].
\end{eqnarray}
For each $\alpha\in (F_{2})^{n}$, we define a linear map 
$$\iota_{\alpha,\x,\y}: A \rightarrow 
\C[[x_{i},x_{i}^{-1},y_{i},y_{i}^{-1}\;|\; i=1,\dots,n]]$$
as follows: First, for $1\le i\le n, \; m\in \Z$, we define
\begin{eqnarray}
\iota_{\alpha,\x,\y}((x_{i}\pm y_{i})^{m})
=\left\{\begin{array}{c}
(x_{i}\pm y_{i})^{m}\;\hspace{3cm}\mbox{
if }\alpha_{i}=0\\
(\pm y_{i}+x_{i})^{m}\hspace{3cm}\mbox{ if }\alpha_{i}=1.
\end{array}\right.
\end{eqnarray}
(We are using the binomial expansion convention.)
Next, we define
\begin{eqnarray}
\iota_{\alpha,\x,\y}((\x\pm \y)^{\m})
=\prod_{i=1}^{n}\iota_{\alpha,\x,\y}((x_{i}\pm y_{i})^{m_{i}})
\end{eqnarray}
for $\m=(m_{1},\dots,m_{n})\in \Z^{n}$. Then we define
\begin{eqnarray}
\iota_{\alpha,\x,\y}\left(\x^{\r}\y^{\s}(\x\pm \y)^{\m}g(\x,\y)\right)
=\x^{\r}\y^{\s}g(\x,\y)\iota_{\alpha,\x,\y}((\x\pm \y)^{\m})
\end{eqnarray}
for $\r,\s,\m\in \Z^{n},\; g(\x,\y)\in {\C}[[\x,\y]]$.
It is clear that each map $\iota_{\alpha,\x,\y}$ is 
linear over the Laurent polynomial algebra in $x_{i},y_{i}$ for $i=1,\dots,n$
and that $\iota_{\alpha,\x,\y}$ commutes with all the formal partial differential 
operators $\partial_{x_{i}}$ and $\partial_{y_{i}}$.

Now, we can write (\ref{eformalproduct}) as:
\begin{eqnarray}\label{edeltajacobin}
\y^{-1}\delta\left(\frac{\x-\z}{\y}\right)
=\sum_{\alpha\in (F_{2})^{n}} (-1)^{|\alpha|}
\iota_{\alpha,\x,\y}\left(\z^{-1}\delta\left(\frac{\x-\y}{\z}\right)\right),
\end{eqnarray}
where $\iota_{\alpha,\x,\y}$ is naturally extended on
 $A[[z_{1},z_{1}^{-1},\dots,z_{n},z_{n}^{-1}]]$.
Furthermore, let $r,s,t\in \Z$ and $p(\x,\y,\z)\in U((\x,\y,\z))$.
Multiplying both sides of (\ref{edeltajacobin}) by $p(\x,\y,\z)$ and then
using the substitution rule (\ref{exzchange}) and (\ref{esubstitution}) we get
\begin{eqnarray}\label{edeltajacobiidentity}
& &\y^{-1}\delta\left(\frac{\x-\z}{\y}\right)p(\y+\z,\y, \z)
\nonumber\\
&=&\sum_{\alpha\in (F_{2})^{n}} (-1)^{|\alpha|}
\iota_{\alpha,\x,\y}\left(\z^{-1}\delta\left(\frac{\x-\y}{\z}\right)
p(\x,\y,\x-\y)\right).
\end{eqnarray}

\br{rvertexstrucureGn}
{\em Set
\begin{eqnarray}
H_{n}=\C[D_{1},\dots, D_{n}],
\end{eqnarray}
where $D_{i}$'s are independent commuting indeterminants or formal variables.
Then $H_{n}$, being identified with the universal enveloping algebra of the abelian
Lie algebra with a basis $\{ D_{1},\dots, D_{n}\}$,
has a cocommutative Hopf algebra structure.
The dual of $H_{n}$ is a commutative associative algebra which 
can be canonically identified as
\begin{eqnarray}
H_{n}^{*}=\C[[x_{1},\dots,x_{n}]].
\end{eqnarray}
Set
\begin{eqnarray}
K_{n}=\C[[x_{1},\dots,x_{n}]][x_{1}^{-1},\dots,x_{n}^{-1}],
\end{eqnarray}
which is a natural commutative associative algebra and a natural 
$H_{n}^{*}$-module. Furthermore, $K_{n}$ is an $H_{n}$-module 
with $D_{i}$ acting as $\partial_{x_{i}}$ for $i=1,\dots,n$.
Define $G_{n}$ to be the pair $(H_{n},K_{n})$ equipped with 
the above structures. From [B2], 
the pair $G_{n}$ is what was called therein an elementary vertex group.
In this paper, we shall at least superficially only use $G_{n}$ as a symbol.}
\er

\bd{dnvertexoperator} {\em (cf. [B2])
Let $W$ be a vector space. A {\em $G_{n}$-vertex operator}
on $W$  is a formal series
$$\psi(\x)=\sum_{\m\in {\Z}^{n}}\psi(\m)\x^{-\m-1}
\in (\End W)[[x_{1},x_{1}^{-1},\dots, x_{n},x_{n}^{-1}]]$$
such that $\psi(\x)w\in W((\x))$ for $w\in W$.}
\ed

All $G_{n}$-vertex operators on $W$ form a subspace, which we
denote by $VO_{G_{n}}(W)$.

\br{rvertexoperatorasamap}
{\em In practice, we shall consider a $G_{n}$-vertex operator 
on $W$ as an object such as a map from $\Z^{n}$ to $\End W$,
which can be uniquely represented by 
a formal series $\psi(\x)$ in {\em any} formal variable $\x$.
In this way, we may use the notations $\psi$, $\psi(\x)$ 
or $\psi(\y)$ for the vertex operator.}
\er

\bd{dlocal1}
{\em $G_{n}$-vertex operators $\psi$ and $\phi$ are said to be
{\em mutually local} if there exists $k\in \N$ such that}
\begin{eqnarray}\label{elocaldef}
(\x-\y)^{k}\psi(\x)\phi(\y)=(\x-\y)^{k}\psi(\y)\psi(\x).
\end{eqnarray}
\ed

A set $S$ of $G_{n}$-vertex operators on $W$ is said to be {\em local}
if any two vertex operators (maybe the same) in $S$ 
are mutually local.

\br{rlocaaldefbigger}
{\em Note that if (\ref{elocaldef}) holds for a certain $k$, 
then it also holds for a bigger $k$.}
\er

\br{rko}
{\em In [KO], a notion of quantum field in one variable and two variables 
were introduced where arbitrary complex powers are allowed.
Note that $VO_{G_{2}}(W)\subset QF_{1}(W)$ and $VO_{G_{4}}(W)\subset QF_{2}(W)$, where
$QF_{1}(W)$ and $QF_{2}(W)$ are the space of quantum fields
in one variable and two variables defined in [KO].}
\er

Let $\psi,\phi\in VO_{G_{n}}(W)$. 
For $w\in W$, $\psi(\x)\phi(\y)w$ involves only 
finitely many negative  powers of all $y_{i}$'s,
but it may involve infinitely many negative powers of $x_{i}$'s.
Therefore, in general, the product
$\psi(\x)\phi(\y)$ does not represent an element of $VO_{G_{2n}}(W)$.

\bl{lproductvertexoperators}
Let $\psi$ and $\phi$  be mutually local $G_{n}$-vertex operators on $W$.
Let $k\in \N$ be such that (\ref{elocaldef}) holds. Then
\begin{eqnarray}\label{eproductvertexoperators}
(\x-\y)^{k}\psi(\x)\phi(\y),\;\; (\x-\y)^{k}\phi(\y)\psi(\x)\in \Hom (W,W((\x,\y))).
\end{eqnarray}
\el

\pf Let $w\in W$. Because
$(\x-\y)^{k}\psi(\x)\phi(\y)w$ involves only finitely many 
negative powers of $y_{i}$'s and
$(\x-\y)^{k}\phi(\y)\psi(\x)w$ involves only finitely many 
negative powers of $x_{i}$'s, by (\ref{elocaldef}) we have
\begin{eqnarray}
(\x-\y)^{k}\psi(\x)\phi(\y)w=(\x-\y)^{k}\phi(\y)\psi(\x)w\in W((\x,\y)).
\end{eqnarray}
This proves (\ref{eproductvertexoperators}). $\;\;\;\;\Box$

\br{rquantumfieldfunc1}
{\em Let $W=\coprod_{\lambda\in \C^{n}}W(\lambda)$ be a $\C^{n}$-graded vector space 
satisfying that for any $\lambda\in \C^{n}$, $W(\m+\lambda)=0$ for 
$\m=(m_{1},\dots,m_{n})\in \Z^{m}$ with $m_{i}$ being
sufficiently small for some $i$. 
Let
$$\psi(\x)=\sum_{\m\in \Z^{n}}\psi(\m)\x^{-\m-1}$$
be a {\em  homogeneous} $G_{n}$-vertex operator
of {\em weight } $h=(h_{1},\dots,h_{n})\in \C^{n}$
in the sense that
\begin{eqnarray}
\psi(\m)W(\lambda)\subset W(\lambda+h-\m-1)\;\;\;\mbox{ for }\m\in \Z^{n},\; 
\lambda\in \C^{n},
\end{eqnarray}
where we consider $1$ as the element $(1,\dots,1)$ of $\C^{n}$.
Let $W'=\coprod_{\lambda\in \C^{n}}W(\lambda)^{*}$ be the restricted dual.
Then
\begin{eqnarray}
\<w', \psi(\x)w\>\in \C[x_{1},x_{1}^{-1},\dots, x_{n},x_{n}^{-1}]
\end{eqnarray}
for $w'\in W',\; w\in W$. Furthermore, let $\psi$ and $\phi$ be mutually 
local homogeneous $G_{n}$-vertex operators on $W$. Then
the formal series 
\begin{eqnarray*}
\<w',\psi(\x)\phi(\y)w\>
\end{eqnarray*}
absolutely converges in the domain $|x_{1}|>|y_{1}|>0,\dots, |x_{n}|>|y_{n}|>0$
to a rational function of the form
\begin{eqnarray}
h(\x,\y)=p(\x,\y)/\left(\x^{\r}\y^{\s}(\x-\y)^{\k}\right)
\end{eqnarray}
where $g(\x,\y)\in {\C}[\x,\y]$ and $\r,\s, \k\in \Z^{n}$.}
\er

Let $R(W,\x,\y)$ consist of formal series
$$f(\x,\y)\in (\End W)[[x_{i},x_{i}^{-1},y_{i},y_{i}^{-1}\;|\; i=1,\dots,n]]$$ 
for which there exists $k\in \N$ such that
\begin{eqnarray}\label{elocaldefhalf}
(\x-\y)^{k}f(\x,\y)\in \Hom (W,W((\x,\y))),
\end{eqnarray}
i.e., $(\x-\y)^{k}f(\x,\y)\in VO_{2n}(W)$.
Clearly, $R(W,\x,\y)$ is a vector subspace.
In view of Lemma \ref{lproductvertexoperators},
if $\psi$ and $\phi$ are mutually local $G_{n}$-vertex operators 
on $W$, then
\begin{eqnarray}
\psi(\x)\phi(\y),\;\;
\phi(\y)\psi(\x)\in R(W,\x,\y).
\end{eqnarray}

\br{rpreparation}
{\em Note that for any $g(\x,\y)\in V_{G_{2n}}(W)$,
\begin{eqnarray}
\iota_{\alpha,\x,\y}\left((\x-\y)^{\m}\right)g(\x,\y)\;\;\;\mbox
{ exists }
\end{eqnarray}
in $(\End W)[[x_{i},x_{i}^{-1},y_{i},y_{i}^{-1}\;|\; i=1,\dots,n]]$ 
for every $\alpha=(\alpha_{1},\dots,\alpha_{n})
\in (F_{2})^{n},\;\m\in {\Z}^{n}$.
Furthermore, for $w\in W$, $\iota_{\alpha,\x,\y}\left((\x-\y)^{\m}\right)g(\x,\y)w$
involves only finitely many
negative powers of $x_{i}$ and $y_{j}$ for $i,j$ with $\alpha_{i}=1$,
$\alpha_{j}=0$ because $\iota_{\alpha,\x,\y}((\x-\y)^{\m})$ has this property.}
\er

\bd{dtoperators}
{\em For $\alpha\in (F_{2})^{n}$, we define a linear map 
$$\iota_{\alpha,\x,\y}: R(W,\x,\y) \rightarrow
(\End W)[[x_{i},x_{i}^{-1}, y_{i},y_{i}^{-1}\;|\;i=1,\dots,n]]$$
by
\begin{eqnarray}\label{edefiota-alpha}
\iota_{\alpha,\x,\y}(f(\x,\y))w
=\iota_{\alpha,\x,\y}\left((\x-\y)^{-k}\right)
\left((\x-\y)^{k}f(\x,\y)w\right)
\end{eqnarray}
for $f(\x,\y)\in R(W,\x,\y),\; w\in W$, 
where $k$ is a nonnegative integer such that (\ref{elocaldefhalf}) holds.}
\ed

In view of Remark \ref{rpreparation}, the expression on the 
right-hand side of (\ref{edefiota-alpha}) exists.
Furthermore, it does not depend on $k$. Indeed, 
let $k_{1}>k_{2}$ be nonnegative integers such that
(\ref{elocaldefhalf}) holds for $k=k_{1},k_{2}$.
Then
\begin{eqnarray*}
& &\iota_{\alpha,\x,\y}\left((\x-\y)^{-k_{1}}\right)
\left((\x-\y)^{k_{1}}f(\x,\y)\right)\nonumber\\
&=&\iota_{\alpha,\x,\y}\left((\x-\y)^{-k_{1}}\right)
(\x-\y)^{k_{1}-k_{2}}\left((\x-\y)^{k_{2}}f(\x,\y)\right)\nonumber\\
&=&\iota_{\alpha,\x,\y}\left((\x-\y)^{-k_{2}}\right)
\left((\x-\y)^{k_{2}}f(\x,\y)\right).
\end{eqnarray*}

\br{rnocancellation}
{\em Note that in general, $\iota_{\alpha,\x,\y}((\x-\y)^{-k})f(\x,\y)$ may not
exist. Thus all the parenthesis in (\ref{edefiota-alpha}) are absolutely necessary.
In general, one should carefully use the associativity to multiply products.}
\er

\br{rexample}
{\em Let $\psi$ and $\phi$ be mutually local $G_{n}$-vertex operators on $W$.
With $\psi(\x)\phi(\y)$ and $\phi(\y)\psi(\x)$ being in $R(W,\x,\y)$,
$\iota_{\alpha,\x,\y}(\psi(\x)\phi(\y))$ 
and $\iota_{\alpha,\x,\y}(\phi(\y)\psi(\x))$ are defined. 
Furthermore, from definition and the locality assumption we have
\begin{eqnarray}\label{e2.32}
\iota_{\alpha,\x,\y}(\psi(\x)\phi(\y))=\iota_{\alpha,\x,\y}(\phi(\y)\psi(\x))
\end{eqnarray}
for every $\alpha\in (F_{2})^{n}$.
Especially, we have
\begin{eqnarray}
\iota_{{\bf 0},\x,\y}(\psi(\x)\phi(\y))&=&\psi(\x)\phi(\y),\\
\iota_{{\bf 1},\x,\y}(\psi(\x)\phi(\y))&=&\phi(\y)\psi(\x),
\end{eqnarray}
where ${\bf 0}=(0,\dots,0)$ and ${\bf 1}=(1,\dots,1)$.}
\er

{}From the definition we immediately have:

\bl{liotamaps}
For $f(\x,\y)\in R(W,\x,\y)$, there exists $k\in \N$ such that
\begin{eqnarray}
(\x-\y)^{k}\iota_{\alpha,\x,\y}(f(\x,\y))=(\x-\y)^{k}f(\x,\y)
\end{eqnarray}
for all $\alpha\in (F_{2})^{n}$. In particular, we may take $k\ge 0$ such that
$$(\x-\y)^{k}f(\x,\y)\in \Hom (W,W((\x,\y))). \;\;\;\;\Box$$
\el

\br{rmoreproperty}
{\em It follows from (\ref{edefiota-alpha}) and Remark \ref{rpreparation}
that for $w\in W$, 
$\iota_{\alpha,\x,\y}(f(\x,\y))w$
involves only finitely many negative (integral) powers of 
$x_{i}$ and $y_{j}$ with $\alpha_{i}=1$ and 
$\alpha_{j}=0$.
For example, $\iota_{e_{i},\x,\y}(\psi(\x)\phi(\y))w$
involves only finitely many negative powers of 
$x_{i}$ and $y_{j}$ for $j\ne i$.}
\er

\bl{lmoreproperty2}
Let $f(\x,\y)\in R(W,\x,\y)$ be such that $(\x-\y)^{\m}f(\x,\y)$ exists
for $\m\in \Z^{n}$. Then
\begin{eqnarray}
\iota_{\alpha,\x,\y}((\x-\y)^{\r+\s}f(\x,\y))
=\iota_{\alpha,\x,\y}((\x-\y)^{\r})\iota_{\alpha,\x,\y}((\x-\y)^{\s}f(\x,\y))
\end{eqnarray}
for $\r,\s\in {\Z}^{n}$. 
\el

\pf In view of Lemma \ref{liotamaps}, there exists $k\in \N$ such that 
\begin{eqnarray}
& &(\x-\y)^{k}\iota_{\alpha,\x,\y}((\x-\y)^{\r})=(\x-\y)^{k+\r},\nonumber\\
& &(\x-\y)^{k}\iota_{\alpha,\x,\y}((\x-\y)^{\s}f(\x,\y))
=(\x-\y)^{k+\s}f(\x,\y),\nonumber\\
& &(\x-\y)^{k}\iota_{\alpha,\x,\y}((\x-\y)^{\r+\s}f(\x,\y))
=(\x-\y)^{\r+\s+k}f(\x,\y).
\end{eqnarray}
Then
\begin{eqnarray}
& &(\x-\y)^{2k}
\iota_{\alpha,\x,\y}((\x-\y)^{\r})\iota_{\alpha,\x,\y}((\x-\y)^{\s}f(\x,\y))
\nonumber\\
&=&(\x-\y)^{k+\r} (\x-\y)^{k+\s}f(\x,\y)\nonumber\\
&=&(\x-\y)^{2k}(\x-\y)^{\r+\s}f(\x,\y)\nonumber\\
&=&(\x-\y)^{2k}\iota_{\alpha,\x,\y}((\x-\y)^{\r+\s}f(\x,\y)).
\end{eqnarray}
Now it follows immediately from Remark \ref{rcancellationfact}.$\;\;\;\;\;\Box$

\br{rquantumfieldfunc2local}
{\em Let $W$ and $W'$ be given as in Remark \ref{rquantumfieldfunc1}. 
Let $\psi, \phi$ be  mutually local homogeneous $G_{n}$-vertex operators.
Let $w'\in W',\; w\in W$ and let $h(\x,\y)$ be as in Remark \ref{rquantumfieldfunc1}. 
Then for each $\alpha\in (F_{2})^{n}$, $\<w',\iota_{\alpha,\x,\y}(\psi(\x)\phi(\y))w\>$
are the formal Laurent series of the same rational function $h(\x,\y)$
in a different domain. }
\er

\br{rtopological}
{\em We here give more information about 
$\iota_{\alpha,\x,\y}(\psi(\x)\phi(\y))$. Recall that for $w\in W$,
\begin{eqnarray}\label{equantitative}
\iota_{\alpha,\x,\y}(\psi(\x)\phi(\y))w=\iota_{\alpha,\x,\y}((\x-\y)^{-k})
\left((\x-\y)^{k}\psi(\x)\phi(\y)w\right),
\end{eqnarray}
where $k$ is a nonnegative integer such that
\begin{eqnarray}
(\x-\y)^{k}\psi(\x)\phi(\y)\in \Hom (W,W((\x,\y))).
\end{eqnarray}
With $(\x-\y)^{k}\psi(\x)\phi(\y)w\in W((\x,\y))$, from (\ref{equantitative})
we see that
each coefficient of formal series $\iota_{\alpha,\x,\y}(\psi(\x)\phi(\y))w$
is a finite sum of the coefficients of $\psi(\x)\phi(\y)w$.
Then the coefficients of $\iota_{\alpha,\x,\y}(\psi(\x)\phi(\y))$ 
lie in the closure of the space spanned by the coefficients of 
$\psi(\x)\phi(\y)$ under a certain natural topology on $\End W$.}
\er

\bl{lcancellationfact}
Let $f(\x,\y)\in R(W,\x,\y)$. Then there exists $k\in \N$ such that
\begin{eqnarray}
(x_{i}-y_{i})^{k}\iota_{\alpha,\x,\y}(f(\x,\y))=
(x_{i}-y_{i})^{k}\iota_{\sigma_{i}(\alpha),\x,\y}(f(\x,\y))
\end{eqnarray}
for $\alpha\in (F_{2})^{n},\; 1\le i\le n.$ 
In particular, we may take $k\ge 0$ such that 
$(\x-\y)^{k}f(\x,\y)\in \Hom (W,W((\x,\y)))$.
\el

\pf Let $k\ge 0$ be such that $(\x-\y)^{k}f(\x,\y)\in \Hom (W,W((\x,\y)))$.
By  Lemma \ref{liotamaps}, we have
\begin{eqnarray}\label{e123}
(\x-\y)^{k}\iota_{\alpha,\x,\y}(f(\x,\y))=(\x-\y)^{k}f(\x,\y)
=(\x-\y)^{k}\iota_{\sigma_{i}(\alpha),\x,\y}(f(\x,\y)).
\end{eqnarray}
Set
\begin{eqnarray}
p_{i}(\x,\y)=\coprod_{j=1, j\ne i}^{n}(x_{j}-y_{j})^{-k}.
\end{eqnarray}
It is easy to see that the associativity assumption holds for both the products
$$\iota_{\alpha,\x,\y}(p_{i}(\x,\y))(\x-\y)^{k}\iota_{\alpha,\x,\y}(f(\x,\y))$$
and
$$\iota_{\alpha,\x,\y}(p_{i}(\x,\y))(\x-\y)^{k}
\iota_{\sigma_{i}(\alpha),\x,\y}(f(\x,\y)).$$
Multiplying both sides of (\ref{e123}) by $\iota_{\alpha,\x,\y}(p_{i}(\x,\y))$ 
and then applying associativity we obtain the desired identity.
$\;\;\;\;\Box$

Furthermore, we have:

\bc{cToperatorproperties}
Let $\psi$ and $\phi$ be mutually local $G_{n}$-vertex operators on $W$.
Then there exists $k\in \N$ such that
\begin{eqnarray}
& &(\x-\y)^{k}\iota_{\alpha,\x,\y}(\psi(\x)\phi(\y))=
(\x-\y)^{k}\psi(\x)\phi(\y),\label{et0010}\\
& &(x_{i}-y_{i})^{k}\iota_{\alpha,\x,\y}(\psi(\x)\phi(\y))=
(x_{i}-y_{i})^{k}\iota_{\sigma_{i}(\alpha),\x,\y}(\psi(\x)\phi(\y))
\label{elocalconni}
\end{eqnarray}
for $\alpha\in (F_{2})^{n},\; 1\le i\le n.$
In particular, we may take $k\ge 0$ such that (\ref{elocaldef}) holds.
\ec

\pf Let $k\in \N$ be such that (\ref{elocaldef}) holds.
Then $$(\x-\y)^{k}\psi(\x)\phi(\y)\in \Hom (W,W((\x,\y))).$$
Now it follows immediately from Lemma \ref{liotamaps} and 
Lemma \ref{lcancellationfact}.$\;\;\;\;\Box$

\bl{llocalproperty1}
Let $\psi$, $\phi$ and $\theta$ be pair-wise mutually 
local $G_{n}$-vertex operators on $W$. 
Then there exists $k\in \N$ such that
\begin{eqnarray}
(\x-\z)^{k}(\y-\z)^{k}\theta (\z)\iota_{\alpha,\x,\y}(\psi(\x)\phi(\y))
=(\x-\z)^{k}(\y-\z)^{k}\iota_{\alpha,\x,\y}(\psi(\x)\phi(\y))\theta (\z)
\end{eqnarray}
for all $\alpha\in (F_{2})^{n}$.
\el

\pf By Corollary \ref{cToperatorproperties} and the locality assumption, 
there exists $k\in \N$ such that (\ref{et0010}) and 
the locality relation  (\ref{elocaldef}) for $(\theta,\psi)$ and $(\theta,\phi)$
hold. Then
\begin{eqnarray}\label{exyxzyz}
& &(\x-\y)^{k}(\x-\z)^{k}(\y-\z)^{k}
\theta(\z)\iota_{\alpha,\x,\y}(\psi(\x)\phi(\y))\nonumber\\
&=&(\x-\y)^{k}(\x-\z)^{k}(\y-\z)^{k}
\theta(\z)\psi(\x)\phi(\y)\nonumber\\
&=&(\x-\y)^{k}(\x-\z)^{k}(\y-\z)^{k}
\psi(\x)\phi(\y)\theta(\z)\nonumber\\
&=&(\x-\y)^{k}(\x-\z)^{k}(\y-\z)^{k}
\iota_{\alpha,\x,\y}(\psi(\x)\phi(\y))\theta(\z).
\end{eqnarray}
Notice that the associativity assumption holds for both the products
$$\iota_{\alpha,\x,\y}((\x-\y)^{-k})(\x-\y)^{k}\left((\x-\z)^{k}(\y-\z)^{k}
\theta(\z)\iota_{\alpha,\x,\y}(\psi(\x)\phi(\y))\right)$$ 
and
$$\iota_{\alpha,\x,\y}((\x-\y)^{-k})(\x-\y)^{k}\left((\x-\z)^{k}(\y-\z)^{k}
\iota_{\alpha,\x,\y}(\psi(\x)\phi(\y))\theta(\z)\right).$$ 
By taking $U=(\End W)[[z_{1},z_{1}^{-1},\dots,z_{n},z_{n}^{-1}]]$ in
Remark \ref{rcancellationfact}, we easily see that we can cancel the factor 
$(\x-\y)^{k}$ from (\ref{exyxzyz}) to get the desired identity.$\;\;\;\;\Box$

\section{Higher dimensional analogues of vertex algebras}

In this section we define a higher dimensional analogue of the notion of
vertex algebra with a Jacobi identity as the main axiom and
we prove basic duality properties. We also give some examples.

\bd{dabstractva}
{\em A {\em $G_{n}$-vertex algebra}
is a vector space $V$ equipped with a linear map 
\begin{eqnarray}
Y: & &V\rightarrow (\End V)[[x_{1},x_{1}^{-1},\dots,x_{n},x_{n}^{-1}]]
\nonumber\\
& &v\mapsto Y(v,\x)=\sum_{\m\in {\Z}^{n}}v_{\m}{\bf x}^{-\m-1}
\;\;\;(\mbox{where }v_{\m}\in \End V)
\end{eqnarray}
and equipped with a distinguished vector ${\bf 1}\in V$ such that 
the following conditions hold:

(A1) For $u\in V$, $Y(u,\x)\in VO_{G_{n}}(V)$, i.e., 
$Y(u,\x)v\in V((\x))$ for $u,v\in V$.

(A2) $Y({\bf 1},\x)=1$ and $Y(v,\x){\bf 1}\in V[[\x]]$ with
$\lim_{\x\rightarrow 0}Y(v,\x){\bf 1}=v$ for $v\in V$.

(A3) For $u,v\in V$, there exists $k\in {\N}$ such that 
\begin{eqnarray}
(\x-\y)^{k}Y(u,\x)Y(v,\y)=(\x-\y)^{k}Y(v,\y)Y(u,\x)
\end{eqnarray}
(the {\em weak commutativity}).

(A4) For $u,v\in V$, 
\begin{eqnarray}\label{ejacobignalgebra}
& &\y^{-1}\delta\left(\frac{\x-\z}{\y}\right)Y(Y(u,\z)v,\y)\nonumber\\
&=&\sum_{\alpha\in (F_{2})^{n}}(-1)^{|\alpha|}\iota_{\alpha,\x,\y}
\left(\z^{-1}\delta\left(\frac{\x-\y}{\z}\right)
Y(u,\x)Y(v,\y)\right)
\end{eqnarray}
(the {\em Jacobi identity}).}
\ed

Note that Axioms (A1) and (A3) are necessary for the Jacobi identity 
(Axiom (A4)) 
to make sense and that $\iota_{\alpha,\x,\y}$ extends naturally on
$R(V,\x,\y)[[z_{1},z_{1}^{-1},\dots,z_{n},z_{n}^{-1}]]$.

\bd{ddoperators}
{\em Let $V$ be a $G_{n}$-vertex algebra. For $1\le i\le n$, we define
$D_{i}\in \End V$ by}
\begin{eqnarray}
D_{i}(v)=\Res_{\x}x_{i}^{-1}\x^{-1}Y(v,\x){\bf 1}\;\;\;\mbox{ for }v\in V.
\end{eqnarray}
\ed

Then from (A2) we have
\begin{eqnarray}\label{eprecreation}
Y(v,\y){\bf 1}\equiv v+y_{1}D_{1}(v)+\cdots +y_{n}D_{n}(v)
\;\;\;\mbox{ mod }\left(\sum_{i=1}^{n}y_{i}^{2} V[[\y]]\right).
\end{eqnarray}
Set
\begin{eqnarray}
\x\D=x_{1}D_{1}+\cdots + x_{n}D_{n}.
\end{eqnarray}

\bp{pbrackerderivative}
Let $V$ be a vertex algebra and let $D_{i}$ be defined as above. Then
\begin{eqnarray}\label{ebrackerderivative}
[D_{i},Y(u,\x)]=Y(D_{i}(u),\x)=\partial_{x_{i}}Y(u,\x)
\end{eqnarray}
for $u\in V$.
\ep

\pf Using the Jacobi identity with $v={\bf 1}$, (A2) and (2.32) we get
\begin{eqnarray}
& &Y(D_{i}(u),\y)\nonumber\\
&=&\Res_{\x}\Res_{\z}z_{i}^{-1}\z^{-1}Y(Y(u,\z){\bf 1},\y)\nonumber\\
&=&\Res_{\x}\Res_{\z}z_{i}^{-1}\z^{-1}
\sum_{\alpha\in (F_{2})^{n}}(-1)^{|\alpha|}
\iota_{\alpha,\x,\y}\left(\z^{-1}\delta\left(\frac{\x-\y}{\z}\right)
Y(u,\x)\right)\nonumber\\
&=&\Res_{\x}\Res_{\z}z_{i}^{-1}\z^{-1}
\left(\sum_{\alpha\in (F_{2})^{n}}(-1)^{|\alpha|}
\iota_{\alpha,\x,\y}\left(\z^{-1}\delta\left(\frac{\x-\y}{\z}\right)\right)
\right)Y(u,\x)\nonumber\\
&=&\Res_{\x}\Res_{\z}z_{i}^{-1}\z^{-1}\y^{-1}
\delta\left(\frac{\x-\z}{\y}\right)
Y(u,\x)\nonumber\\
&=&\Res_{\x}\Res_{\z}z_{i}^{-1}\z^{-1}
\x^{-1}\delta\left(\frac{\y+\z}{\x}\right)
Y(u,\x)\nonumber\\
&=&\Res_{\z}z_{i}^{-1}\z^{-1}Y(u,\y+\z)\nonumber\\
&=&\Res_{\z}z_{i}^{-1}\z^{-1}e^{(z_{1}\partial_{y_{1}}+\cdots 
+z_{n}\partial_{y_{n}})}Y(u,\y)\nonumber\\
&=&\partial_{y_{i}}Y(u,\y).
\end{eqnarray}
This proves the second equality of (\ref{ebrackerderivative}).

Applying the Jacobi identity to ${\bf 1}$ and then taking $\Res_{\x}$ 
and using (\ref{eprecreation}) we get
\begin{eqnarray}
& &D_{i}(Y(u,\z)v)\nonumber\\
&=&\Res_{\y}y_{i}^{-1}\y^{-1}Y(Y(u,\z)v,\y){\bf 1}\nonumber\\
&=&\Res_{\x}\Res_{\y}y_{i}^{-1}\y^{-1}\sum_{\alpha\in (F_{2})^{n}} (-1)^{|\alpha|}
\iota_{\alpha,\x,\y}\left(\z^{-1}\delta\left(\frac{\x-\y}{\z}\right)\right)
\iota_{\alpha,\x,\y}(Y(u,\x)Y(v,\y){\bf 1})\nonumber\\
&=&\Res_{\x}\Res_{\y}y_{i}^{-1}\y^{-1}\sum_{\alpha\in (F_{2})^{n}} (-1)^{|\alpha|}
\iota_{\alpha,\x,\y}\left(\z^{-1}\delta\left(\frac{\x-\y}{\z}\right)\right)
\iota_{\alpha,\x,\y}(Y(v,\y)Y(u,\x){\bf 1})\nonumber\\
&=&\Res_{\x}\Res_{\y}y_{i}^{-1}\y^{-1}\z^{-1}\delta\left(\frac{\x-\y}{\z}\right)
Y(u,\x)Y(v,\y){\bf 1}\nonumber\\
&=&\Res_{\y}y_{i}^{-1}\y^{-1}Y(u,\z+\y)Y(v,\y){\bf 1}\nonumber\\
&=&\Res_{\y}y_{i}^{-1}\y^{-1}
e^{(y_{1}\partial_{z_{1}}+\cdots +y_{n}\partial_{z_{n}})}
Y(u,\z)Y(v,\y){\bf 1}\nonumber\\
&=&Y(u,\z)D_{i}(v)+\partial_{z_{i}}Y(u,\z)v,
\end{eqnarray}
noting that because for $m\in \Z$, $\iota_{\alpha,\x,\y}(\x-\y)^{m}$
involves only nonnegative powers of $x_{i}$ for $i$ 
with $\alpha_{i}=1$,
$$\iota_{\alpha,\x,\y}((\x-\y)^{m})
\iota_{\alpha,\x,\y}(Y(v,\y)Y(u,\x){\bf 1})$$
involves only nonnegative powers of $x_{i}$ for $i$ with $\alpha_{i}=1$, 
\begin{eqnarray}
& &\Res_{\x}\iota_{\alpha,\x,\y}
\left((\x-\y)^{m}Y(u,\x)Y(v,\y){\bf 1}\right)
\nonumber\\
&=&\Res_{\x}\iota_{\alpha,\x,\y}((\x-\y)^{m})
\iota_{\alpha,\x,\y}(Y(u,\x)Y(v,\y){\bf 1})\nonumber\\
&=&\Res_{\x}\iota_{\alpha,\x,\y}((\x-\y)^{m})
\iota_{\alpha,\x,\y}(Y(v,\y)Y(u,\x){\bf 1})\nonumber\\
&=&0
\end{eqnarray}
for all $m\in \Z,\; \alpha\ne (0,\dots,0)$. Then
\begin{eqnarray}
D_{i}Y(u,\z)v-Y(u,\z)D_{i}(v)=\partial_{z_{i}}Y(u,\z)v.
\end{eqnarray}
This proves the outside equality of (\ref{ebrackerderivative}).
Now the proof is complete.$\;\;\;\;\Box$

{}From (\ref{ebrackerderivative}) and the Taylor Theorem we have 
the following conjugation formula:
\begin{eqnarray}\label{edconjugation}
e^{\y \D}Y(u,\x)e^{-\y\D}=Y(u,\x+\y)\;\;\;\mbox{ for }u\in V.
\end{eqnarray}

We also have the following skew-symmetry:

\bp{pskew-symmetry}
In the setting of Proposition \ref{pbrackerderivative} we have
\begin{eqnarray}\label{eskewsymmetry}
Y(u,\x)v=e^{\x\D}Y(v,-\x)u\;\;\;\mbox{ for }u,v\in V.
\end{eqnarray}
\ep

\pf Define a translation  $\sigma$ of $(F_{2})^{n}$ by
\begin{eqnarray}
\sigma(\alpha)=\alpha+(1,\dots,1)\;\;\;\mbox{ for }\alpha\in (F_{2})^{n}.
\end{eqnarray}
Then $(-1)^{|\sigma(\alpha)|}=(-1)^{n}(-1)^{|\alpha|}$ and
\begin{eqnarray}
\iota_{\alpha,\y,\x}=\iota_{\sigma(\alpha),\x,\y}\;\;\;
\mbox{ for }\alpha\in (F_{2})^{n}.
\end{eqnarray}
It follows that the right-hand side of the Jacobi identity stays the same
if we replace $(\x,\y,\z)$ by $(\y,\x,-\z)$.
Just like in the ordinary case [FHL], using this, the Jacobi identity 
and the creation property we obtain the skew-symmetry (\ref{eskewsymmetry}).
$\;\;\;\;\Box$

Similar to the ordinary case, we have the following weak associativity:

\bp{pweakcommassoc}
Let $V$ be a vertex algebra. For $u,w\in V$, there exists $l\ge 0$ such that
\begin{eqnarray}\label{eweakassociativityequiv}
(\z+\y)^{l}Y(u,\z+\y)Y(v,\y)w=
(\z+\y)^{l} Y(Y(u,\z)v,\y)w
\end{eqnarray}
for all $v\in V$. In particular, we may take $l\ge 0$ such that
$\x^{l}Y(u,\x)w\in V[[\x]]$. On the other hand, in 
Definition \ref{dabstractva}
let us only assume (A1)-(A3), in addition we assume the weak associativity property.
Then the Jacobi identity holds, so that $V$ is a vertex algebra.
\ep

\pf Assume the Jacobi identity. For any $l\in \N$,
applying $\Res_{\x} \x^{l}$ to both sides of the Jacobi identity 
and using delta function substitution property we get
\begin{eqnarray}
& &(\y+\z)^{l}Y(Y(u,\z)v,\y)\nonumber\\
&=&\sum_{\alpha\in (F_{2})^{n}}(-1)^{|\alpha|}\Res_{\x}\iota_{\alpha,\x,\y}
\left( \z^{-1}\delta\left(\frac{\x-\y}{\z}\right)\x^{l}Y(u,\x)Y(v,\y)w\right).
\end{eqnarray}
Let $l\in \N$ be such that $\x^{l}Y(u,\x)w\in V[[\x]]$.
Recall (\ref{e2.32}):
$$\iota_{\alpha,\x,\y}(Y(u,\x)Y(v,\y))=
\iota_{\alpha,\x,\y}(Y(v,\y)Y(u,\x)).$$
Because for $m\in \Z$, $\iota_{\alpha,\x,\y}(\x-\y)^{m}$
involves only nonnegative powers of $x_{i}$ for $i$ 
with $\alpha_{i}=1$,
$$\iota_{\alpha,\x,\y}((\x-\y)^{m})
\iota_{\alpha,\x,\y}(\x^{l}Y(v,\y)Y(u,\x)w)$$
involves only nonnegative powers of $x_{i}$ for $i$ with $\alpha_{i}=1$. 
Then
\begin{eqnarray}
& &\Res_{\x}\iota_{\alpha,\x,\y}
\left((\x-\y)^{m}\x^{l}Y(u,\x)Y(v,\y)w\right)
\nonumber\\
&=&\Res_{\x}\iota_{\alpha,\x,\y}((\x-\y)^{m})
\iota_{\alpha,\x,\y}(\x^{l}Y(u,\x)Y(v,\y)w)\nonumber\\
&=&\Res_{\x}\iota_{\alpha,\x,\y}((\x-\y)^{m})
\iota_{\alpha,\x,\y}(\x^{l}Y(v,\y)Y(u,\x)w)\nonumber\\
&=&0
\end{eqnarray}
for all $m\in \Z,\; \alpha\ne (0,\dots,0)$.
Therefore
\begin{eqnarray}
& &(\y+\z)^{l}(Y(Y(u,\z)v,\y)w\nonumber\\
&=&\Res_{\x}
\z^{-1}\delta\left(\frac{\x-\y}{\z}\right)
\x^{l}Y(u,\x)Y(v,\y)w\nonumber\\
&=&\Res_{\x}
\x^{-1}\delta\left(\frac{\z+\y}{\x}\right)
(\z+\y)^{l}Y(u,\x)Y(v,\y)w\nonumber\\
&=&(\z+\y)^{l}Y(u,\z+\y)Y(v,\y)w.
\end{eqnarray}

Now assume (\ref{eweakassociativityequiv}). 
Set
\begin{eqnarray}
p(\x,\y)=\x^{l}(\x-\y)^{k}Y(u,\x)Y(v,\y)w\in W((\x,\y))
\end{eqnarray}
and
\begin{eqnarray}
h(\x,\y,\z)=\x^{-l}\z^{-k}p(\x,\y)\in W((\x,\y,\z)).
\end{eqnarray}
Then
\begin{eqnarray}
& &Y(u,\x)Y(v,\y)w=h(\x,\y,\x-\y),\\
& &Y(Y(u,\z)v,\y)w=h(\y+\z,\y,z).
\end{eqnarray}
Then the Jacobi identity immediately follows from (\ref{edeltajacobiidentity}).
$\;\;\;\;\Box$
 
Taking $\Res_{\z}$ of the Jacobi identity we get
\begin{eqnarray}\label{ecommutatorformulanot}
\sum_{\alpha\in (F_{2})^{n}}(-1)^{|\alpha|}\iota_{\alpha,\x,\y}(Y(u,\x)Y(v,\y))
=\Res_{\z}\y^{-1}\delta\left(\frac{\x-\z}{\y}\right)Y(Y(u,\z)v,\y).
\end{eqnarray}
This is a generalization of Borcherds' commutator formula.
By taking $\Res_{\x}$ we also have the following iterate formula:
\begin{eqnarray}\label{eiterate}
Y(Y(u,\z)v,\y)=\sum_{\alpha\in (F_{2})^{n}}(-1)^{|\alpha|}
\Res_{\x}\iota_{\alpha,\x,\y}
\left( \z^{-1}\delta\left(\frac{\x-\y}{\z}\right)Y(u,\x)Y(v,\y)\right).
\end{eqnarray}
Note that the above iterate formula gives rise to an associator formula:
\begin{eqnarray}\label{eassociatorformula}
& &Y(Y(u,\z)v,\y)-Y(u,\z+\y)Y(v,\y)\nonumber\\
&=&\sum_{\alpha\in (F_{2})^{n},\alpha\ne {\bf 0}}(-1)^{|\alpha|}
\Res_{\x}\iota_{\alpha,\x,\y}
\left( \z^{-1}\delta\left(\frac{\x-\y}{\z}\right)Y(u,\x)Y(v,\y)\right).
\end{eqnarray}
But (\ref{ecommutatorformulanot})
does not give a genuine commutator formula.
Nevertheless, similar to the ordinary vertex algebra case, we have:

\bp{pcommproperty}
Let $u,v\in V$. Then $[Y(u,\x),Y(v,\y)]=0$ if and only if
$Y(u,\x)v\in V[[\x]]$.
\ep

\pf Suppose that $[Y(u,\x),Y(v,\y)]=0$. Then with $Y(u,\x){\bf 1}\in V[[\x]]$, 
we have
$$Y(u,\x)v=\Res_{\y}\y^{-1}Y(u,\x)Y(v,\y){\bf 1}=\Res_{\y}\y^{-1}
Y(v,\y)Y(u,\x){\bf 1}\in V[[\x]].$$

Conversely, assume $Y(u,\x)v\in V[[\x]]$. 
In view of Proposition \ref{pweakcommassoc}, for any $w\in V$, we have
\begin{eqnarray}\label{eweakassocl=0}
Y(u,\z+\y)Y(w,\y)v=Y(Y(u,\z)w,\y)v.
\end{eqnarray}
Applying $e^{-\y \D}$ from left to (\ref{eweakassocl=0}) and then using 
the conjugation formula (\ref{edconjugation}) and 
the skew-symmetry (\ref{eskewsymmetry}) (twice) we obtain
\begin{eqnarray}
Y(u,\z)Y(v,-\y)w=Y(v,-\y)Y(u,\z)w.
\end{eqnarray}
This proves that $[Y(u,\x),Y(v,\y)]=0$.
Now the proof is complete.$\;\;\;\;\Box$

Similar to the ordinary case ([FHL], [DL], [Li2], cf. [Ka]) (in fact, using 
exactly the same argument of [Li2]) we have:

\bt{tequivalentaxioms}
Let $(V,Y,{\bf 1},\D)$ be a quadruple consisting of a vector space $V$,
a vector ${\bf 1}$ in $V$, a linear map $Y$ from $V$ to
$(\End V)[[x_{1},x_{1}^{-1},\dots,x_{n},x_{n}^{-1}]]$ and a vector 
$$\D=(D_{1},\dots,D_{n})\in (\End V)^{n}$$
such that $D_{i}({\bf 1})=0$ for $i=1,\dots,n$ and such that
the following conditions hold:

(E1) $Y(u,\x)v\in V((\x))$ for $u,v\in V$.

(E2) $Y({\bf 1},\x)=1$ and $Y(v,\x){\bf 1}\in V[[\x]]$ with 
$\lim_{\x\mapsto 0}Y(v,\x){\bf 1}=v$ for $v\in V$.

(E3) $[D_{i},Y(u,\x)]=\partial_{x_{i}}Y(u,\x)$.

(E4) For $u,v\in V$, there exists $k\in \N$ such that
$$(\x-\y)^{k}Y(u,\x)Y(v,\y)=(\x-\y)^{k}Y(v,\y)Y(u,\x).$$
Then $(V,Y,{\bf 1})$ is a vertex algebra and 
$D_{i}(v)=\Res_{\x}x_{i}^{-1}\x^{-1}Y(v,\x){\bf 1}$ for $v\in V$.
\et

\pf First, from (E3) and the Taylor Theorem we immediately have
(cf. (\ref{edconjugation}))
\begin{eqnarray}\label{eDconjugate}
e^{\y \D}Y(v,\x)e^{-\y \D}=Y(v,\x+\y)\;\;\;\mbox{ for }v\in V.
\end{eqnarray}
Applying (\ref{eDconjugate}) to ${\bf 1}$, using the fact that 
$D_{i}({\bf 1})= 0$ for $i=1,\dots,n$, we get 
\begin{eqnarray}
e^{\y \D}Y(v,\x){\bf 1}=Y(v,\x+\y){\bf 1}.
\end{eqnarray}
Since $Y(v,\x+\y){\bf 1}\in V[[(\x+\y)]]$, we may set $\y=-\x$. 
Then using the creation property we get
$e^{-\x \D}Y(v,\x){\bf 1}=v$. That is,
\begin{eqnarray}
Y(v,\x){\bf 1}=e^{\x \D}v.
\end{eqnarray}
Second, we have the skew-symmetry as follows:
Let $k\in \N$ be such that the weak commutativity holds and such that
$\x^{k}Y(v,\x)u\in V[[\x]]$. Then
\begin{eqnarray}
(\x-\y)^{k}Y(u,\x)Y(v,\y){\bf 1}
&=&(\x-\y)^{k}Y(v,\y)Y(u,\x){\bf 1}\nonumber\\
&=&(\x-\y)^{k}Y(v,\y)e^{\x \D}u\nonumber\\
&=&(\x-\y)^{k}e^{\x \D}Y(v,\y-\x)u.
\end{eqnarray}
Since $(\x-\y)^{k}Y(v,\y-\x)u$ involves only nonnegative powers 
of $(\x-\y)$, we may
set $\y$ to zero. Setting $\y=0$ and using creation property we get
\begin{eqnarray}
\x^{k}Y(u,\x)v=\x^{k}e^{\x \D}Y(v,-\x)u,
\end{eqnarray}
which immediately gives (cf. (\ref{eskewsymmetry}))
\begin{eqnarray}
Y(u,\x)v=e^{\x \D}Y(v,-\x)u\;\;\;\mbox{ for }u,v\in V.
\end{eqnarray}
Third, let $l\in \N$ be such that the weak commutativity holds 
for $(u,w)$. Then using the skew-symmetry and conjugation formula 
we obtain the weak associativity as
\begin{eqnarray}
(\z+\y)^{l}Y(u,\z+\y)Y(v,\y)w&=&(\z+\y)^{l}Y(u,\z+\y)e^{-\y \D}Y(w,-\y)v
\nonumber\\
&=&(\z+\y)^{l}e^{-\y \D}Y(u,\z)Y(w,-\y)v\nonumber\\
&=&(\z+\y)^{l}e^{-\y \D}Y(w,-\y)Y(u,\z)v\nonumber\\
&=&(\z+\y)^{l}Y(Y(u,\z)v,\y)w.
\end{eqnarray}
Now it follows from the second part of Proposition \ref{pweakcommassoc}
that $(V,Y,{\bf 1})$ is a vertex algebra. The last assertion follows from
(E3) and the assumption that $D_{i}({\bf 1})=0$ for $i=1,\dots,n$.
$\;\;\;\;\Box$

\br{rjacobinotcommassoc}
{\em As proved in [DL], [FHL] and [Li2] for the ordinary case, 
in view of Theorem \ref{tequivalentaxioms}, one may
define the notion of vertex algebra by using Conditions (E1)-(E4) 
as the axioms (see for example [Ka] in the ordinary case).}
\er

\br{rdifferent}
{\em In fact, in Theorem \ref{tequivalentaxioms},
(E1) follows from (E2) and (E4), so that (E1) is redundant.
Indeed, for $u,v\in V$, let $k\in \N$ be such that the weak 
commutativity relation holds for $(u,v)$. Then
\begin{eqnarray}
(\x-\y)^{k}Y(u,\x)Y(v,\y){\bf 1}=(\x-\y)^{k}Y(v,\y)Y(u,\x){\bf 1}.
\end{eqnarray}
Since the expression on the left-hand side involves only 
nonnegative powers of $y_{i}$'s and 
the expression on the right-hand side involves only 
nonnegative powers of $x_{i}$'s, we have
\begin{eqnarray}
(\x-\y)^{k}Y(u,\x)Y(v,\y){\bf 1}\in V[[\x,\y]].
\end{eqnarray}
Setting $\y=0$, then using creation property, we get
\begin{eqnarray}
\x^{k}Y(u,\x)v\in V[[\x]].
\end{eqnarray}
That is, $Y(u,\x)v\in V((\x))$.}
\er

\bp{pskewassoctocomm}
In the definition of the notion $G_{n}$-vertex algebra, let us only
assume (A1)-(A2) and in addition we assume that the skew-symmetry
(\ref{eskewsymmetry}) and weak associativity hold. 
Then weak commutativity and the Jacobi
identity hold.
\ep

\pf Let $u,v,w\in W$ and let $l\in \N$ be such that
\begin{eqnarray}
(\z+\y)^{l}Y(u,\z+\y)Y(v,\y)w=(\z+\y)^{l}Y(Y(u,\z)v,y)w.
\end{eqnarray}
Applying $e^{-\y\D}$ from left and then using 
skew-symmetry (\ref{eskewsymmetry}) we get
\begin{eqnarray}\label{eassoctocommproof}
& &(\z+\y)^{l}e^{-\y\D}Y(u,\z+\y)e^{\y\D}Y(w,-\y)v\nonumber\\
&=&(\z+\y)^{l}Y(w,-\y)Y(u,\z)v.
\end{eqnarray}
In particular, for $w={\bf 1}$, we obtain
\begin{eqnarray}
(\z+\y)^{l}e^{-\y\D}Y(u,\z+\y)e^{\y\D}v=(\z+\y)^{l}Y(u,\z)v,
\end{eqnarray}
which immediately implies the conjugation formula (\ref{edconjugation}).

Now, we go back to let $w$ be arbitrary again. From (\ref{eassoctocommproof})
using the conjugation formula (\ref{edconjugation}) we get
\begin{eqnarray}
(\z+\y)^{l}Y(u,\z)Y(w,-\y)v=(\z+\y)^{l}Y(w,-\y)Y(u,\z)v.
\end{eqnarray}
Notice that $l$ depends only on $u$ and $w$, not $v$. This gives 
weak commutativity. The rest of the assertions immediately follows.
$\;\;\;\;\Box$

\br{rrelations}
{\em Any $G_{n}$-vertex algebra $V$ is a natural $G_{n+1}$-vertex algebra.
Indeed, Axioms (A1)-(A3) clearly hold and for $u,v\in V$. Furthermore, 
multiplying 
both sides of (\ref{ejacobignalgebra}) by 
$y_{n+1}^{-1}\delta\left(\frac{x_{n+1}-z_{n+1}}{y_{n+1}}\right)$ 
and then using the identity 
\begin{eqnarray}
y_{n+1}^{-1}\delta\left(\frac{x_{n+1}-z_{n+1}}{y_{n+1}}\right)
=z_{n+1}^{-1}\delta\left(\frac{x_{n+1}-y_{n+1}}{z_{n+1}}\right)
-z_{n+1}^{-1}\delta\left(\frac{-y_{n+1}+x_{n+1}}{z_{n+1}}\right)
\end{eqnarray}
to rewrite the new right-hand side we obtain the Jacobi identity
for $G_{n+1}$-vertex algebras.

On the other hand, let $(V,Y,{\bf 1})$ be a $G_{n}$-vertex algebra. 
For any subset $I$ of $\{1,\dots,n\}$ we set
\begin{eqnarray}
V_{I}=\{ v\in V\;|\; D_{i}(v)=0\;\;\;\mbox{ for }i\in I\}.
\end{eqnarray}
In view of Proposition \ref{pbrackerderivative}, for $v\in V$, $v\in V_{I}$ 
if and only if 
\begin{eqnarray}
\partial_{x_{i}}Y(v,\x)=0\;\;\;\mbox{ for }i\in I.
\end{eqnarray}
Notice that for $u,v\in V_{I}$,
\begin{eqnarray}
\left(\coprod_{i\in I}\Res_{x_{i}}\right) \iota_{\alpha,\x,\y}\left(
\z^{-1}\delta\left(\frac{\x-\y}{\z}\right)Y(u,\x)Y(v,\y)\right)=0
\end{eqnarray}
for $\alpha\in (F_{2})^{n}$ with $\alpha_{i}=1$ for some $i\in I$.
Then by taking $\coprod_{i\in I}\Res_{x_{i}}$ of (\ref{ejacobignalgebra})
we obtain the Jacobi identity for $G_{n-|I|}$-vertex algebras.
Then it follows that $(V_{I},Y,{\bf 1})$ is a $G_{n-|I|}$-vertex algebra.}
\er

\bd{dmodule}
{\em Let $V$ be a vertex algebra. A $V$-{\em module} is a vector space $W$ 
equipped with a linear map $Y$ from $V$ to 
$(\End W)[[x_{1},x^{-1},\dots, x_{n},x_{n}^{-1}]]$
that satisfies the following conditions:

(M1) $Y(v,\x)w\in W((\x))$ for $v\in V,\; w\in W$.

(M2) $Y({\bf 1},\x)=1$.

(M3) The weak commutativity (A3) on $W$ holds.

(M4) The Jacobi identity (\ref{ejacobignalgebra}) on $W$ holds.}
\ed

Recall that Jacobi identity implies weak associativity.
In fact, in the notion of module, weak associativity
implies weak commutativity and Jacobi identity (cf. [Le], [Li2]).

\bp{pmodduleassociativity}
Let $V$ be a vertex algebra and let $W$ be a vector space equipped with 
a linear map $Y$ from $V$ to 
$(\End W)[[x_{1},x_{1}^{-1},\dots,x_{n},x_{n}^{-1}]]$
such that (M1) and (M2) of Definition \ref{dmodule} hold and such that
the following weak associativity holds: For every $u\in V,\; w\in W$,
there exists $l \in \N$ such that
\begin{eqnarray}
(\z+\y)^{l}Y(u,\z+\y)Y(v,\y)w=(\z+\y)^{l}Y(Y(u,\z)v,y)w
\end{eqnarray}
for all $v\in V$. Then $W$ is a $V$-module.
\ep

\pf In view of Proposition \ref{pweakcommassoc} (the second part),
it suffices to prove the weak commutativity.

Let $u,v\in V,\; w\in W$. From (M1) and the given weak
associativity, there exists
$k\in \N$ (only depending on $u,v$) such that 
$$\z^{k}Y(u,\z)v\in V[[\z]]$$
and there exists
$l\in \N$ such that all the following conditions hold:
\begin{eqnarray}
& & \y^{l}Y(v,\y)w\in W[[\y]],\\
& &(\z+\y)^{l}Y(u,\z+\y)Y(v,\y)w=(\z+\y)^{l}Y(Y(u,\z)v,\y)w,\\
& &(\x-\z)^{l}Y(v,-\z+\x)Y(u,\x)w=(\x-\z)^{l}Y(Y(v,-\z)u,\x)w.
\end{eqnarray}
Then 
$$\y^{l}\z^{k}(\z+\y)^{l}Y(u,\z+\y)Y(v,\y)w
=\y^{l}\z^{k}(\z+\y)^{l}Y(Y(u,\z)v,\y)w\in W[[\y,\z]].$$
Set
\begin{eqnarray}
p(\y,\z)=\y^{l}\z^{k}(\z+\y)^{l}Y(u,\z+\y)Y(v,\y)w.
\end{eqnarray}
Then
\begin{eqnarray}
\y^{l}(\x-\y)^{k}\x^{l}Y(u,\x)Y(v,\y)w=p(\y,\x-\y)
\end{eqnarray}
because
\begin{eqnarray}
& &\z^{-1}\delta\left(\frac{\x-\y}{\z}\right)
\x^{l}\y^{l}(\x-\y)^{k}Y(u,\x)Y(v,\y)w\nonumber\\
&=&\x^{-1}\delta\left(\frac{\z+\y}{\x}\right)
\y^{l}\z^{k}(\z+\y)^{l}Y(u,\z+\y)Y(v,\y)w\nonumber\\
&=&\x^{-1}\delta\left(\frac{\z+\y}{\x}\right)p(\y,\z)\nonumber\\
&=&\x^{-1}\delta\left(\frac{\x-\y}{\z}\right)p(\y,\x-\y).
\end{eqnarray}
On the other hand, using skew-symmetry and Taylor theorem we obtain
\begin{eqnarray}
& &\z^{-1}\delta\left(\frac{\y-\x}{-\z}\right)
\x^{l}\y^{l}(\x-\y)^{k}Y(v,\y)Y(u,\x)w\nonumber\\
&=&\y^{-1}\delta\left(\frac{-\z+\x}{\y}\right)
\x^{l}(-\z+\x)^{l}\z^{k}Y(v,-\z+\x)Y(u,\x)w\nonumber\\
&=&\y^{-1}\delta\left(\frac{-\z+\x}{\y}\right)
\left(\x^{l}(-\z+\x)^{l}\z^{k}Y(Y(v,-\z)u,\x)w\right)\nonumber\\
&=&\y^{-1}\delta\left(\frac{-\z+\x}{\y}\right)
\left(\x^{l}\z^{k}(-\z+\x)^{l}Y(e^{-\z \D}Y(u,\z)v,\x)w\right)\nonumber\\
&=&\y^{-1}\delta\left(\frac{-\z+\x}{\y}\right)
\left(\x^{l}\z^{k}(-\z+\x)^{l}Y(Y(u,\z)v,\x-\z)w\right)\nonumber\\
&=&\y^{-1}\delta\left(\frac{-\z+\x}{\y}\right)
\left(\y^{l}\z^{l}(\y+\z)^{l}Y(Y(u,\z)v,\y)w\right)|_{\y=\x-\z}\nonumber\\
&=&\y^{-1}\delta\left(\frac{-\z+\x}{\y}\right)p(\y,\z)|_{\y=\x-\z}.
\end{eqnarray}
Since $p(\y,\z)\in W[[\y,\z]]$, 
$$p(\y,\z)|_{\y=\x-\z}=p(\y,\z)|_{\y=-\z+\x}.$$
Then
\begin{eqnarray}
& &\z^{-1}\delta\left(\frac{\y-\x}{-\z}\right)
\x^{l}\y^{l}(\x-\y)^{l}Y(v,\y)Y(u,\x)w\nonumber\\
&=&\y^{-1}\delta\left(\frac{-\z+\x}{\y}\right)
p(\y,\z)|_{\y=-\z+\x}\nonumber\\
&=&\y^{-1}\delta\left(\frac{-\z+\x}{\y}\right)p(\y,\z)\nonumber\\
&=&\z^{-1}\delta\left(\frac{\y-\x}{-\z}\right)p(\y,-\y+\x)\nonumber\\
&=&\z^{-1}\delta\left(\frac{\y-\x}{-\z}\right)p(\y,\x-\y).
\end{eqnarray}
Then
\begin{eqnarray}
\x^{l}\y^{l}(\x-\y)^{k}Y(v,\y)Y(u,\x)w=p(\y,-\y+\x).
\end{eqnarray}
Consequently,
\begin{eqnarray}
\x^{l}\y^{l}(\x-\y)^{k}Y(u,\x)Y(v,\y)w=\x^{l}\y^{l}(\x-\y)^{k}Y(v,\y)Y(u,\x)w.
\end{eqnarray}
Thus
$$(\x-\y)^{k}Y(u,\x)Y(v,\y)w=(\x-\y)^{k}Y(v,\y)Y(u,\x)w.$$
Recall that $k$ does not depend on $w$. This
proves the weak commutativity and concludes the proof.$\;\;\;\;\Box$

\br{rmodulecomments}
{\em Given a $G_{n}$-vertex algebra $V$, we have seen that
in the notion of $V$-module, weak associativity is stronger
than weak commutativity. This is contrary to the situation
for the notion of $G_{n}$-vertex algebra where
the vacuum vector plays a crucial role in the proof.}
\er

We later shall need the following result (cf. [Li1]):

\bl{lvacuumlikeweak}
Let $V$ be a $G_{n}$-vertex algebra and $S$ a generating set of $V$.
Let $W$ be a $V$-module and let $w\in W$ be such that
$Y(u,\x)w\in W[[\x]]$ for $u\in S$. Then
$Y(v,\x)w\in V[[\x]]$ for all $v\in V$.
Furthermore, the linear map 
\begin{eqnarray}
f: V\rightarrow W;\; v\mapsto  \lim_{\x\rightarrow 0}Y(v,\x)w
\end{eqnarray}
is a $V$-homomorphism.
\el

\pf For the first part, it suffices to prove that all elements $v$ of $V$
satisfying $Y(v,\x)w\in W[[\x]]$ form a vertex subalgebra of $V$.
Furthermore, it is enough to prove that if $Y(u,\x)w,\; Y(v,\x)w\in W[[\x]]$,
then $Y(u_{\m}v,\x)w\in W[[\x]]$ for $\m\in \Z^{n}$.
In view of Proposition \ref{pweakcommassoc} we have
\begin{eqnarray}\label{e3.58}
Y(Y(u,\z)v,\y)w=Y(u,\z+\y)Y(v,\y)w\;\;\;\mbox{ for }u,v\in V.
\end{eqnarray}
Then $Y(u_{\m}v,\x)w\in W[[\x]]$ for $\m\in \Z^{n}$.

By taking $\lim_{\x\rightarrow 0}$ in (\ref{e3.58}) we obtain
\begin{eqnarray}
f(Y(u,\z)v)=Y(u,\z)f(v)\;\;\;\mbox{ for }u,v\in V.
\end{eqnarray}
That is, $f$ is a $V$-homomorphism from $V$ to $W$.$\;\;\;\;\Box$


Furthermore we have:

\bp{pisomorphism}
Let $V$, $S$, $W$ and $w$ be given as in Lemma \ref{lvacuumlikeweak}.
Furthermore, assume that $W$ is a faithful $V$-module in the sense
that $Y$ is one-to-one, $w$ generates $W$ as a $V$-module 
and assume that there exists 
$\d=(d_{1},\dots,d_{n})\in (\End W)^{n}$ such that
$d_{i}(w)=0$ and
\begin{eqnarray}\label{edbracketisom}
[d_{i},Y(v,\x)]=\partial_{x_{i}}Y(v,\x)
\end{eqnarray}
on $W$ for $i=1,\dots,n, \; v\in V$. Then the linear map $f$ defined
in Lemma \ref{lvacuumlikeweak} is a $V$-isomorphism.
\ep

\pf In view of Lemma \ref{lvacuumlikeweak}, from the assumption that 
$w$ generates $W$ as a $V$-module, $f$ is an onto $V$-homomorphism.
It remains to prove that $f$ is one-to-one. 
Assume that $f(v)=0$ for some $v\in V$. From the assumptions we get
\begin{eqnarray}
Y(v,\x+\y)w=e^{\y\D}Y(v,\x)e^{-\y\D}w=e^{\y\D}Y(v,\x)w.
\end{eqnarray}
With $Y(v,\x+\y)w\in W[[(\x+\y)]]$, we may set $\x=\0$ to obtain  
$Y(v,\y)w=e^{\y\D}f(v)=0$.
Let $a\in W$ be such that $Y(v,\x)a=0$ and let $u\in S$. 
Then there exists $k\in \N$ such that
$$(\x-\y)^{k}Y(v,\x)Y(u,\y)a=(\x-\y)^{k}Y(u,\y)Y(v,\x)a=0.$$
Thus $Y(v,\x)Y(u,\y)a=0$.
Then it follows from induction that 
$Y(v,\x)=0$ on $W$. With $W$ being faithful, we must have $v=0$.
This proves that $f$ is one-to-one.
$\;\;\;\;\Box$

{\em Note:} In fact, in Proposition \ref{pisomorphism}, it is enough to
assume that (\ref{edbracketisom}) holds for $v\in S$. (See the proof
Proposition \ref{paddderivation}.)

In the following we shall give some examples of 
$G_{n}$-vertex algebras. 

\bexa{ecommassocder}
{\em Let $A$ be a commutative associative algebra with identity $1$ equipped 
with $n$ (possibly the same) pairwise commuting derivations 
$D_{1},\dots, D_{n}$. 
Set
\begin{eqnarray}
{\bf D}=(D_{1},\dots,D_{n}).
\end{eqnarray}
Then as a convention we write
\begin{eqnarray}
\x {\bf D}=x_{1}D_{1}+\cdots +x_{n}D_{n}.
\end{eqnarray}
For $a\in A$, we define 
$$Y(a,\x)\in (\End A)[[x_{1},\dots,x_{n}]]
\subset (\End A)[[x_{1},x_{1}^{-1},\dots,x_{n},x_{n}^{-1}]]$$
by
\begin{eqnarray}
Y(a,\x)b=\left(e^{\x {\bf D}}a\right)b=
\left(e^{x_{1}D_{1}+\cdots+x_{n}D_{n}}a\right)b
\end{eqnarray}
for $b\in A$. Then
\begin{eqnarray}	
Y(a,\x)Y(b,\y)c
=\left(e^{x_{1}D_{1}+\cdots+x_{n}D_{n}}a\right)
\left(e^{y_{1}D_{1}+\cdots+y_{n}D_{n}}b\right)c
\end{eqnarray}
for $a,b,c\in A$.
Since $D_{1},\dots,D_{n}$ are pairwise commuting, we have
\begin{eqnarray}	
Y(a,\x)Y(b,\y)=Y(b,\y)Y(a,\x)\;\;\;\mbox{ for }a,b\in A.
\end{eqnarray}
In particular, the weak commutativity holds.
Clearly, the vacuum properties hold with ${\bf 1}=1$.
Then
\begin{eqnarray}
\iota_{\alpha,\x,\y}(Y(a,\x)Y(b,\y))
&=&Y(a,\x)Y(b,\y)\nonumber\\
&=&\left(e^{x_{1}D_{1}+\cdots+x_{n}D_{n}}a\right)
\left(e^{y_{1}D_{1}+\cdots+y_{n}D_{n}}b\right).
\end{eqnarray}
On the other hand,
\begin{eqnarray}
Y(Y(a,\z)b,\y)c
&=&\left(e^{y_{1}D_{1}+\cdots+y_{n}D_{n}}
\left(e^{z_{1}D_{1}+\cdots+z_{n}D_{n}}a\right)b
\right)c\nonumber\\
&=&\left(e^{(y_{1}+z_{1})D_{1}+\cdots+(y_{n}+z_{n})D_{n}}a\right)
\left(e^{y_{1}D_{1}+\cdots+y_{n}D_{n}}b\right)c.
\end{eqnarray}
Then 
\begin{eqnarray}
& &\y^{-1}\delta\left(\frac{\x-\z}{\y}\right)
Y(Y(a,\z)b,\y)c\nonumber\\
&=&\y^{-1}\delta\left(\frac{\x-\z}{\y}\right)
\left(e^{(\y+\z){\bf D}}a\right)
\left(e^{\y{\bf D}}b\right)c\nonumber\\
&=&\x^{-1}\delta\left(\frac{\y+\z}{\x}\right)
\left(e^{\x{\bf D}}a\right)
\left(e^{\y{\bf D}}b\right)c.
\end{eqnarray}
Then the Jacobi identity immediately follows.}
\eexa

\br{rholo}
{\em Let $V$ be a $G_{n}$-vertex algebra such that
\begin{eqnarray}
Y(v,\x)\in (\End V)[[x_{1},\dots,x_{n}]]\;\;\;\mbox{ for all }v\in V.
\end{eqnarray}
In view of Proposition \ref{pcommproperty} we have
\begin{eqnarray}
[Y(u,\x),Y(v,\y)]=0\;\;\;\mbox{ for }u,v\in V.
\end{eqnarray}
In view of Proposition \ref{pweakcommassoc} we have
\begin{eqnarray}
Y(Y(u,\z)v,\y)w=Y(u,\z+\y)Y(v,\y)w\;\;\;\mbox{ for }u,v,w\in V.
\end{eqnarray}
Then one can easily see that $V$ equipped with the product
$$u\cdot v=\Res_{\x}\x^{-1}Y(u,\x)v\;\;\;\mbox{ for }u,v\in V$$
is a commutative
associative algebra with ${\bf 1}$ as identity. Furthermore, 
$D_{i}$ are mutually commuting derivations where
$$D_{i}(v)=\Res_{\x}\x^{-1}x_{i}^{-1}Y(v,\x){\bf 1}\;\;\;\mbox{ for }v\in V.$$
Then $Y(u,\x)v=\left(e^{\x {\bf D}}u\right)v$ for $u,v\in V$.}
\er

\bexa{exatensorproduct}
{\em Let $V_{1},\dots,V_{n}$ be ordinary vertex algebras. Set
\begin{eqnarray}
V=V_{1}\otimes \cdots \otimes V_{n}
\end{eqnarray}
as a vector space. Define
\begin{eqnarray}
Y(v,\x)=Y(v^{1},x_{1})\otimes \cdots\otimes Y(v^{n},x_{n})
\in (\End V)[[x_{1},x_{1}^{-1},\dots,x_{n},x_{n}^{-1}]]
\end{eqnarray}
for $v=v^{1}\otimes \cdots\otimes v^{n}\in V$.
Use linearity to extend $Y$ to be a linear map from  $V$ to 
$(\End V)[[x_{1},x_{1}^{-1},\dots,x_{n},x_{n}^{-1}]]$.
Set
\begin{eqnarray}
{\bf 1}={\bf 1}\otimes \cdots\otimes {\bf 1}\in V.
\end{eqnarray}
Now we show that $(V, Y, {\bf 1})$ is a $G_{n}$-vertex algebra.
Clearly, the lower truncation condition and the vacuum property hold.

For $u=u^{1}\otimes \cdots\otimes u^{n}, \;
v=v^{1}\otimes \cdots\otimes v^{n}\in V$, we have
\begin{eqnarray}
Y(u,\x)Y(v,\y)
=Y(u^{1},x_{1})Y(v^{1},y_{1})\otimes \cdots \otimes 
Y(u^{n},x_{n})Y(v^{n},y_{n}).
\end{eqnarray}
Then weak commutativity obviously holds.
Furthermore, we have
\begin{eqnarray}
& &Y(Y((u^{1}\otimes\cdots\otimes u^{n}),\z)(v^{1}\otimes\cdots\otimes v^{n}),\y)
\nonumber\\
&=&Y(Y(u^{1},z_{1})v^{1}\otimes\cdots\otimes Y(u^{n},z_{n})v^{n},\y)\nonumber\\
&=&Y(Y(u^{1},z_{1})v^{1},y_{1})\otimes\cdots\otimes Y(Y(u^{n},z_{n})v^{n},y_{n}).
\end{eqnarray}

For two formal series $A(x), B(x)\in (\End U)[[x,x^{-1}]]$, we define
\begin{eqnarray}
\tau_{0}(A(x),B(y))=A(x)B(y),\;\;\;\tau_{1}(A(x),B(y))=B(y)A(x).
\end{eqnarray}
Then for $\a=(\a_{1},\dots,\a_{n})\in (F_{2})^{n}$,
\begin{eqnarray}
\iota_{\a,\x,\y}(Y(u,\x)Y(v,\y))=
\tau_{\a_{1}}(Y(u^{1},x_{1}),Y(v^{1},y_{1}))\otimes \cdots \otimes 
\tau_{\a_{n}}(Y(u^{n},x_{n}),Y(v^{n},y_{n})).
\end{eqnarray}
With this, we easily see that the Jacobi identity is exactly 
the tensor product of $n$ ordinary Jacobi identities.

Furthermore, let $W_{i}$ be a $V_{i}$-module for $i=1,\dots,n$. Then
$W_{1}\otimes\cdots \otimes W_{n}$ is a $V$-module with the obvious action.}
\eexa

By considering certain subalgebras of the
tensor product $G_{n}$-vertex algebras
we can obtain $G_{n}$-vertex algebras a little bit more general than the
tensor product $G_{n}$-vertex algebras.

\bexa{exafixedpoint}
{\em Let $G$ be an abelian group and let
$U=\oplus_{g\in G}U_{g}$ and $V=\oplus_{g\in G}V_{g}$
ordinary vertex algebras graded by $G$ such that
\begin{eqnarray}
& &u_{\m}u'\in U_{g+h}\;\;\;\mbox{ for }u\in U_{g},\; u'\in U_{h},
\; g,h\in G,\\
& &v_{\m}v'\in V_{g+h} \;\;\;\mbox{ for }v\in V_{g},
\; v'\in V_{h},\; g,h\in G.
\end{eqnarray}
Set
\begin{eqnarray}
(U\otimes V)_{G}=\oplus_{g\in G} U_{g}\otimes V_{g}.
\end{eqnarray}
It is easy to see that $(U\otimes V)_{G}$ is a vertex subalgebra 
of $U\otimes V$.}
\eexa

A {\em homomorphism} from a $G_{n}$-vertex algebra $U$ to another $V$ is a
linear map $f$ such that $f({\bf 1})={\bf 1}$ and 
\begin{eqnarray}
f(Y(u,\x)v)=Y(f(u),\x)f(v)\;\;\;\mbox{ for }u,v\in U.
\end{eqnarray}
Isomorphism and automorphism are defined in the obvious ways.

Just as in the case for an ordinary vertex algebra, from
a routine proof we immediately have:

\bl{lfixedpoint}
Let $V$ be a $G_{n}$-vertex algebra and 
$G$ a group of automorphisms of $V$. Set
\begin{eqnarray}
V^{G}=\{ v\in V\;|\; g(v)=v\;\;\;\mbox{ for }g\in G\}.
\end{eqnarray}
Then $V^{G}$ is a vertex subalgebra. $\;\;\;\;\Box$
\el

\section{Vertex algebras of $G_{n}$-vertex operators}

In this section, we shall continue the study of Section 2
on $G_{n}$-vertex operators on a vector space $W$.
For any two (maybe the same) given mutually local $G_{n}$-vertex operators,
we define a family of products parametrized by $\m\in \Z^{n}$.
As our main results we prove that any maximal local space of 
$G_{n}$-vertex operators on $W$ equipped with the defined products
is a vertex algebra  with $W$ as a natural module
and that any set of pair-wise mutually local $G_{n}$-vertex operators 
on $W$ automatically generates a vertex algebra.

Now, we continue the discussion of Section 2 and
we shall freely use the notations and conventions defined there.
Let $W$ be a vector space. In terms of maps $\iota_{\alpha,\x,\y}$, 
we introduce the following definition.

\bd{dvertexmult}
{\em Let $\psi$ and $\phi$ be mutually local $G_{n}$-vertex operators on $W$. 
For $\m\in {\Z}^{n}$, we define a formal series 
$$(\psi_{\m}\phi)(\y)\in (\End W)[[y_{1},y_{1}^{-1},\dots,y_{n},y_{n}^{-1}]]$$
 by
\begin{eqnarray}
& &(\psi_{\m}\phi)(\y)\nonumber\\
&=&\sum_{\alpha\in (F_{2})^{n}} (-1)^{|\alpha|}\Res_{\x}\Res_{\z}
\iota_{\alpha,\x,\y}\left(\z^{-1}\delta\left(\frac{\x-\y}{\z}\right)
\z^{\m}\psi(\x)\phi(\y)\right)\nonumber\\
&=&\sum_{\alpha\in (F_{2})^{n}}(-1)^{|\alpha|}
\Res_{\x}\iota_{\alpha,\x,\y}\left((\x-\y)^{\m}
\psi(\x)\phi(\y)\right).\label{edefcomponents1}
\end{eqnarray}}
\ed

Noting that if $\m\ge 0$, for every $\alpha\in (F_{2})^{n}$, 
$$\iota_{\alpha,\x,\y}((\x-\y)^{\m})=(\x-\y)^{\m},$$
we get
\begin{eqnarray}
(\psi_{\m}\phi)(\y)
&=&\Res_{\x}(\x-\y)^{\m}
\sum_{\alpha\in (F_{2})^{n}}(-1)^{|\alpha|}
\iota_{\alpha,\x,\y}(\psi(\x)\phi(\y)).
\end{eqnarray}
(The sum is analogous to the usual commutator.)

Let $\alpha=(\alpha_{1},\dots,\alpha_{n})\in (F_{2})^{n},\; w\in W$.
Consider $\Res_{\x}\iota_{\alpha,\x,\y}\left((\x-\y)^{\m}
\psi(\x)\phi(\y)w\right)$.
Let $k\in \N$ be such that (\ref{elocaldef}) holds. Then
\begin{eqnarray}
& &\iota_{\alpha,\x,\y}\left((\x-\y)^{\m}\psi(\x)\phi(\y)w\right)\nonumber\\
&=&\iota_{\alpha,\x,\y}((\x-\y)^{\m-k})
\left((\x-\y)^{k}\psi(\x)\phi(\y)w\right).
\end{eqnarray}
For $1\le i\le n$, if $\alpha_{i}=0$, both $(\x-\y)^{k}\psi(\x)\phi(\y)w$
and $\iota_{\alpha,\x,\y}((\x-\y)^{\m-k})$ involves only finitely many negative 
powers of $y_{i}$ (recall Remark \ref{rmoreproperty}), 
so does $\Res_{\x}\iota_{\alpha,\x,\y}\left((\x-\y)^{\m}
\psi(\x)\phi(\y)w\right)$.
If $\alpha_{i}=1$, we have 
$$\iota_{\alpha,\x,\y}((\x-\y)^{\m-k})
=(-y_{i}+x_{i})^{m_{i}-k}\prod_{j\ne i}\iota_{\alpha,\x,\y}((x_{j}-y_{j})^{m_{j}-k}).$$
Then 
$$\Res_{\x}\iota_{\alpha,\x,\y}((\x-\y)^{\m-k})
\left((\x-\y)^{k}\psi(\x)\phi(\y)w\right)$$
involves only finitely many negative powers of $y_{i}$, so does 
$\Res_{\x}\iota_{\alpha,\x,\y}\left((\x-\y)^{\m}
\psi(\x)\phi(\y)w\right)$. Thus, $\Res_{\x}\iota_{\alpha,\x,\y}\left((\x-\y)^{\m}
\psi(\x)\phi(\y)w\right)$ always involves only finitely many 
negative powers of $y_{i}$'s. From (\ref{edefcomponents1}) 
we have
$$(\psi_{\m}\phi)(\y)w\in W((\y))\;\;\;\mbox{ for }w\in W.$$
Thus,
$(\psi_{\m}\phi)(\y)\in VO_{G_{n}}(W)$.
Assume  $m_{i}\ge k$ for some $i$.
Using Corollary \ref{cToperatorproperties} we get
\begin{eqnarray}
& &\iota_{\alpha,\x,\y}\left((\x-\y)^{\m}\psi(\x)\phi(\y)\right)\nonumber\\
&=&\left(\prod_{j\ne i}\iota_{\alpha,\x,\y}((x_{j}-y_{j})^{m_{j}})\right)
(x_{i}-y_{i})^{m_{i}}
\iota_{\alpha,\x,\y}(\psi(\x)\phi(\y))\nonumber\\
&=&\left(\prod_{j\ne i}\iota_{\alpha,\x,\y}((x_{j}-y_{j})^{m_{j}})\right)
(x_{i}-y_{i})^{m_{i}}
\iota_{\sigma_{i}(\alpha),\x,\y}(\psi(\x)\phi(\y))\nonumber\\
&=&\iota_{\sigma_{i}(\alpha),\x,\y}\left((\x-\y)^{\m}\psi(\x)\phi(\y)\right).
\end{eqnarray}
Then using the relation $(-1)^{|\alpha|}=-(-1)^{|\sigma_{i}(\alpha)|}$, we get
$\psi_{\m}\phi=0$. To summarize, we have:

\bl{lbaiscproperty1}
Let $\psi$ and $\phi$ be mutually local $G_{n}$-vertex operators on $W$. Then
\begin{eqnarray}
\psi_{\m}\phi\in VO_{G_{n}}(W)\;\;\;\mbox{  for }\m\in {\Z}^{n}. 
\end{eqnarray}
Furthermore,
\begin{eqnarray}
\psi_{\m}\phi=0\;\;\;\mbox{ if }m_{i}\ge k\;\;\mbox{ for some } i,
\end{eqnarray}
where $\m=(m_{1},\dots,m_{n})\in \Z^{n}$ and $k$ is a nonnegative integer 
such that (\ref{elocaldef}) holds.
$\;\;\;\;\Box$
\el

Set
\begin{eqnarray}\label{egeneratingfunctiondef}
Y(\psi,\z)\phi =\sum_{\m\in \Z^{n}}(\psi_{\m}\phi) \z^{\m-1}.
\end{eqnarray}
In terms of this generating function we have
\begin{eqnarray}\label{eassociator1}
(Y(\psi,\z)\phi)(\y)
=\sum_{\alpha\in (F_{2})^{n}}(-1)^{|\alpha|}\Res_{\x}
\iota_{\alpha,\x,\y}\left(\z^{-1}\delta\left(\frac{\x-\y}{\z}\right)
\psi(\x)\phi(\y)\right).
\end{eqnarray}

\bl{ldefweakassociativity}
Let $\psi,\phi$ be mutually local $G_{n}$-vertex operators on $W$ and let
$w\in W$. Then there exists $l\in \N$ (depending only on $\psi$ and $w$)
such that
\begin{eqnarray}\label{eweakassociativityabstract}
(\y+\z)^{l}(Y(\psi,\z)\phi)(\y)w=(\z+\y)^{l}\psi(\z+\y)\phi(\y)w.
\end{eqnarray}
In particular, we may take $l\ge 0$ such that $\x^{l}\psi(\x)w\in W[[[\x]]$.
\el

\pf For any $l\in \N$, using (\ref{eassociator1}) and 
delta function substitution property we have
\begin{eqnarray}
& &(\y+\z)^{l}(Y(\psi,\z)\phi)(\y)w\nonumber\\
&=&\sum_{\alpha\in (F_{2})^{n}}(-1)^{|\alpha|}\Res_{\x}
(\y+\z)^{l}\iota_{\alpha,\x,\y}\left(\z^{-1}\delta\left(\frac{\x-\y}{\z}\right)
\psi(\x)\phi(\y)w\right)\nonumber\\
&=&\sum_{\alpha\in (F_{2})^{n}}(-1)^{|\alpha|}\Res_{\x}
\iota_{\alpha,\x,\y}\left(\z^{-1}\delta\left(\frac{\x-\y}{\z}\right)
\x^{l}\psi(\x)\phi(\y)w\right).
\end{eqnarray}
Let $l\ge 0$ be such that $\x^{l}\psi(\x)w\in W[[[\x]]$.
Recall (\ref{e2.32}):
$$\iota_{\alpha,\x,\y}(\psi(\x)\phi(\y))=
\iota_{\alpha,\x,\y}(\phi(\y)\psi(\x)).$$
Because for $m\in \Z$, $\iota_{\alpha,\x,\y}(\x-\y)^{m}$
involves only nonnegative powers of $x_{i}$ for $i$ 
with $\alpha_{i}=1$,
$$\iota_{\alpha,\x,\y}((\x-\y)^{\m})\iota_{\alpha,\x,\y}(\x^{l}\phi(\y)\psi(\x)w)$$
involves only nonnegative powers of $x_{i}$ for $i$ with $\alpha_{i}=1$. Then
\begin{eqnarray}
& &\Res_{\x}\iota_{\alpha,\x,\y}\left((\x-\y)^{m}\x^{l}\psi(\x)\phi(\y)w\right)
\nonumber\\
&=&\Res_{\x}\iota_{\alpha,\x,\y}((\x-\y)^{m})
\iota_{\alpha,\x,\y}(\x^{l}\psi(\x)\phi(\y)w)\nonumber\\
&=&\Res_{\x}\iota_{\alpha,\x,\y}((\x-\y)^{m})
\iota_{\alpha,\x,\y}(\x^{l}\phi(\y)\psi(\x)w)\nonumber\\
&=&0
\end{eqnarray}
for all $m\in \Z,\; \alpha\ne (0,\dots,0)$.
Therefore
\begin{eqnarray}
& &(\y+\z)^{l}(Y(\psi,\z)\phi)(\y)w\nonumber\\
&=&\Res_{\x}\x^{l}
\z^{-1}\delta\left(\frac{\x-\y}{\z}\right)
\psi(\x)\phi(\y)w\nonumber\\
&=&\Res_{\x}(\z+\y)^{l}
\x^{-1}\delta\left(\frac{\z+\y}{\x}\right)
\psi(\x)\phi(\y)w\nonumber\\
&=&(\z+\y)^{l}\psi(\z+\y)\phi(\y)w.
\end{eqnarray}
This completes the proof.
$\;\;\;\;\Box$

By observing both sides of (\ref{eweakassociativityabstract}) we 
see that the common quantity on the both sides lies in
$W((\y,\z))$, so that
$$\psi(\z+\y)\phi(\y)w\in W((\y,\z))[(\z+\y)^{-1}].$$
Then
\begin{eqnarray}
(\y+\z)^{l}(Y(\psi,\z)\phi)(\y)w
&=&\iota_{\0, \y,\z}((\z+\y)^{l}\psi(\z+\y)\phi(\y)w)\nonumber\\
&=&(\y+\z)^{l}\iota_{\0, \y,\z}(\psi(\z+\y)\phi(\y)w).
\end{eqnarray}
Thus
\begin{eqnarray}
(Y(\psi,\z)\phi)(\y)w=\iota_{\0, \y,\z}(\psi(\z+\y)\phi(\y)w).
\end{eqnarray}
Therefore, we have proved:

\bc{cassocdefprod}
Let $\psi$ and $\phi$ be mutually local $G_{n}$-vertex operators on $W$. Then
\begin{eqnarray}\label{eassocdefprod}
(Y(\psi,\z)\phi)(\y)=\iota_{\0, \y,\z}(\psi(\z+\y)\phi(\y)).\;\;\;\;\Box
\end{eqnarray}
\ec

\br{rassocdefprod}
{\em We may use equation (\ref{eassocdefprod}) alternatively to define
$\psi_{\m}\phi$.}
\er

With Lemma \ref{ldefweakassociativity}, from
Proposition \ref{pweakcommassoc} we immediately have:

\bp{pjacobirep}
Let $\psi,\phi$ be mutually local $G_{n}$-vertex operators on $W$. Then
\begin{eqnarray}
& &\y^{-1}\delta\left(\frac{\x-\z}{\y}\right)(Y(\psi,\z)\phi)(\y)\nonumber\\
&=&\sum_{\alpha\in (F_{2})^{n}}(-1)^{|\alpha|}
\iota_{\alpha,\x,\y}\left(\z^{-1}\delta\left(\frac{\x-\y}{\z}\right)
\psi(\x)\phi(\y)\right).\;\;\;\;\Box
\end{eqnarray}
\ep

We shall need the following result:

\bl{llocalityprelemma}
Let $\psi,\phi$ and $\theta$ be pairwise mutually local $G_{n}$-vertex 
operators on $W$. Then there exists a nonnegative integer $k$ such that
\begin{eqnarray}\label{eprelemma}
(\x-\y-\z)^{k}(\x-\y)^{k}\theta(\x)(Y(\psi,\z)\phi)(\y)
=(\x-\y-\z)^{k}(\x-\y)^{k}(Y(\psi,\z)\phi)(\y)\theta(\x).
\end{eqnarray}
In particular, we may take $k\ge 0$ such that
\begin{eqnarray}\label{ekconditions}
(\x-\y)^{k}[\theta(\x),\psi(\y)]=0,\;\; (\x-\y)^{k}[\theta(\x),\phi(\y)]=0.
\end{eqnarray}
\el

\pf Let $k\ge 0$ be such that (\ref{ekconditions}) holds.
Using Corollary \ref{cassocdefprod}, 
noticing that $\iota_{\0,\y,\z}$ is linear on $\C[\x,\y,\z]$ we obtain
(\ref{eprelemma}) as
\begin{eqnarray}
& &(\x-\y-\z)^{k}(\x-\y)^{k}\theta(\x)(Y(\psi,\z)\phi)(\y)\nonumber\\
&=&\iota_{\0, \y,\z}\left((\x-\y-\z)^{k}(\x-\y)^{k}
\theta(\x)\psi(\z+\y)\phi(\y)\right)
\nonumber\\
&=&\iota_{\0, \y,\z}\left((\x-\y-\z)^{k}(\x-\y)^{k}\psi(\z+\y)\phi(\y)\theta(\x)
\right)\nonumber\\
&=&(\x-\y-\z)^{k}(\x-\y)^{k}\iota_{\0, \y,\z}(\psi(\z+\y)\phi(\y))\theta(\x)
\nonumber\\
&=&(\x-\y-\z)^{k}(\x-\y)^{k}(Y(\psi,\z)\phi)(\y)\theta(\x).\;\;\;\;\;\Box
\end{eqnarray}

The following is the key result of this section (cf. [Li2]):

\bp{plocaalitygenerating}
Let $\psi$, $\phi$ and $\theta$ be pairwise mutually local 
$G_{n}$-vertex operators on $W$. Then for any $\m\in {\Z}^{n}$,
$\psi_{\m}\phi$ and $\theta$ are mutually local. 
\ep

\pf Let $k\ge 0$ be such that (\ref{eprelemma}) holds.
For any fixed $\m=(m_{1},\dots,m_{n})\in \Z^{n}$, 
in view of Lemma \ref{lbaiscproperty1}, there exists $r\ge 0$ such that
$z_{i}^{r+m_{i}}Y(\psi,\z)\phi$ involves only nonnegative powers of $z_{i}$
for $i=1,\dots,n$, so that
\begin{eqnarray}
\Res_{\z}(\x-\y-\z)^{k+r-\i}
\z^{\m+\i}(\x-\y)^{k}(Y(\psi,\z)\phi)(\y)=0
\end{eqnarray}
for $\i=i_{1},\dots,i_{n})\in \N^{n}$ with $i_{j}\ge r$ for some $j$.
Then using Lemma \ref{llocalityprelemma} we obtain
\begin{eqnarray}
& &(\x-\y)^{2k+r}\theta(\x)(\psi_{\m}\phi)(\y)\nonumber\\
&=&\Res_{\z}\z^{\m}(\x-\y)^{2k+r}\theta(\x)(Y(\psi,\z)\phi)(\y)\nonumber\\
&=&\Res_{\z}\sum_{\i\in \N^{n}}{k+r\choose \i}(\x-\y-\z)^{k+r-\i}
\z^{\m+\i}(\x-\y)^{k}\theta(\x)(Y(\psi,\z)\phi)(\y)\nonumber\\
&=&\Res_{\z}\sum_{\i\in \N^{n}, \i\le r}{k+r\choose \i}(\x-\y-\z)^{k+r-\i}
\z^{\m+\i}(\x-\y)^{k}\theta(\x)(Y(\psi,\z)\phi)(\y)\nonumber\\
&=&\Res_{\z}\sum_{\i\in \N^{n}, \i\le r}{k+r\choose \i}(\x-\y-\z)^{k+r-\i}
\z^{\m+\i}(\x-\y)^{k}(Y(\psi,\z)\phi)(\y)\theta(\x)\nonumber\\
&=&\Res_{\z}\sum_{\i\in \N^{n}}{k+r\choose \i}(\x-\y-\z)^{k+r-\i}
\z^{\m+\i}(\x-\y)^{k}(Y(\psi,\z)\phi)(\y)\theta(\x)\nonumber\\
&=&\Res_{\z}\z^{\m}(\x-\y)^{2k+r}(Y(\psi,\z)\phi)(\y)\theta(\x)\nonumber\\
&=&(\x-\y)^{2k+r}(\psi_{\m}\phi)(\y)\theta(\x).
\end{eqnarray}
This proves that $\psi_{\m}\phi$ and $\theta$ are mutually local.
$\;\;\;\;\Box$

\br{roldcase}
{\em Motivated by the ordinary case [Li2], 
for $\m\in \Z^{n}$, we consider the formal series defined by
\begin{eqnarray}
(\psi_{[\m]}\phi)(\y)
=\Res_{\x}\left((\x-\y)^{\m}\psi(\x)\phi(\y)
-(-\y+\x)^{\m}\phi(\y)\psi(\x)\right).
\end{eqnarray}
Clearly, $\psi_{[\m]}\phi\in VO_{G_{n}}(W)$.
However, in general, $\psi_{[\m]}\phi$ will not be local with the old
vertex operators $\psi$ and $\phi$. Notice that
Corollary \ref{cToperatorproperties} is crucial in the proof of 
Proposition \ref{plocaalitygenerating}.}
\er

\bp{paddderivation}
Let $W$ be a vector space and $\d=(d_{1},\dots,d_{n})$ an $n$-tuple
of endomorphisms of $W$.
Let $VO_{G_{n}}(W,\d)$ consist of $G_{n}$-vertex operators $\psi$ on $W$ such that
\begin{eqnarray}
[d_{i},\psi(\x)]=\partial_{x_{i}}\psi(\x)\;\;\;\mbox{ for }i=1,\dots,n.
\end{eqnarray}
If $\psi,\phi\in VO_{G_{n}}(W,\d)$ are mutually local, then 
\begin{eqnarray}\label{elocalsystemwithd}
\psi_{\m}\phi\in VO_{G_{n}}(W,\d)\;\;\;\mbox{ for }\m\in \Z^{n}.
\end{eqnarray}
\ep

\pf Using Corollary \ref{cassocdefprod} we obtain
\begin{eqnarray}
& &[d_{i},(Y(\psi,\z)\phi)(\y)]\nonumber\\
&=&[d_{i},\iota_{\0,\y,\z}(\psi(\z+\y)\phi(\y))]\nonumber\\
&=&\iota_{\0,\y,\z}\left([d_{i},\psi(\z+\y)]\phi(\y)+\psi(\z+\y)[d_{i},\phi(\y)]\right)
\nonumber\\
&=&\iota_{\0,\y,\z}\left(\left(\partial_{y_{i}}\psi(\z+\y)\right)
\phi(\y)+\psi(\z+\y)\partial_{y_{i}}\phi(\y)\right)
\nonumber\\
&=&\partial_{y_{i}}\iota_{\0,\y,\z}(\psi(\z+\y)\phi(\y))\nonumber\\
&=&\partial_{y_{i}}(Y(\psi,\z)\phi)(\y).
\end{eqnarray}
In terms of components we have
\begin{eqnarray}
[d_{i},(\psi_{\m}\phi)(\y)]=\partial_{y_{i}}(\psi_{\m}\phi)(\y)
\end{eqnarray}
for $i=1,\dots,n$ and for $\m\in \Z^{n}$.
This proves  (\ref{elocalsystemwithd}).$\;\;\;\;\Box$

There are $n$ linear operators on $VO_{G_{n}}(W)$:
\begin{eqnarray}
D_{i}=\partial_{x_{i}}\;\left(={\partial\over\partial x_{i}}\right)
\;\;\; \mbox{ for }i=1,\dots,n.
\end{eqnarray}
Set
\begin{eqnarray}
\D =(D_{1},\dots,D_{n})\in (\End VO_{G_{n}}(W))^{n}.
\end{eqnarray}
Write
\begin{eqnarray}
\x\D=x_{1}D_{1}+\cdots +x_{n}D_{n}.
\end{eqnarray}
Note that every vertex operator on $W$ is mutually local with 
the identity operator $1$ of $W$.

\bl{ldoperatorvacuum}
Let $\psi\in VO_{G_{n}}(W)$. Then
\begin{eqnarray}\label{eabstractvaccum}	
(Y(\psi,\z)1)(\y)=e^{\z \D}\psi(\y).
\end{eqnarray}
In particular, for $1\le i\le n$,
\begin{eqnarray}
D_{i}(\psi)=\Res_{\x}\x^{-1}x_{i}^{-1}Y(\psi,\x)1.
\end{eqnarray}
\el

\pf Using Corollary \ref{cassocdefprod} we obtain (\ref{eabstractvaccum}) as	
\begin{eqnarray}
(Y(\psi,\z)1)(\y)
&=&\iota_{\0,\y,\z}(\psi(\z+\y)1(\y))\nonumber\\
&=&\iota_{\0,\y,\z}(\psi(\z+\y))\nonumber\\
&=&\psi(\y+\z)\nonumber\\
&=&e^{\z \D}\psi(\y).\;\;\;\;\Box
\end{eqnarray}


\bd{dclosedalgebra}
{\em A local subspace $S$ of $VO_{G_{n}}(W)$ is said to be {\em closed} if 
\begin{eqnarray}
\phi_{\m}\phi\in S\;\;\;\mbox{ for }\psi,\phi\in S,\;\m\in {\Z}^{n}.
\end{eqnarray}
By the term ``{\em a vertex algebra of mutually local 
$G_{n}$-vertex operators on $W$}''
we mean a closed local subspace of $VO_{G_{n}}(W)$ that contains 
the identity operator $1$.}
\ed

In the following we shall justify the use of ``vertex algebra'' in the above definition
by proving that a vertex algebra of mutually local $G_{n}$-vertex operators on $W$
is a vertex algebra as defined in Section 3.

Since
\begin{eqnarray}
\Res_{\x}\x^{-1}x_{i}^{-1}(Y(\psi,\x)1)(\y)
=\partial_{x_{i}}\psi(\x)
\end{eqnarray}
for $i=1,\dots,n$, any closed local subspace of $VO_{G_{n}}(W)$
is stable under the action of $D_{i}=\partial_{x_{i}}$ 
for $i=1,\dots,n$.

As an immediate consequence of Proposition \ref{plocaalitygenerating} we have:

\bc{cmaximalclosed}
Any maximal local subspace (or subset) of $VO_{G_{n}}(W)$ is closed. $\;\;\;\;\Box$
\ec

\br{rmaaximal}
{\em Let $S$ be any local subset of $VO_{G_{n}}(W)$. Consider all
the local subsets that contain $S$. It follows from Zorn's lemma that
there exists a maximal local subspace of $VO_{G_{n}}(W)$ that contains $S$.
In view of Corollary \ref{cmaximalclosed}, there exists a closed
local subspace of $VO_{G_{n}}(W)$ that contains $S$ and $1$.
Let $\<S\>$ be the smallest closed local subspace of $VO_{G_{n}}(W)$
that contains $S$ and the identity operator $1$ of $W$. Then $\<S\>$ is 
closed.}
\er

Furthermore, we have:

\bp{pdoperatorproperty}
Let $V$ be a closed local subspace of $VO_{G_{n}}(W)$. Then
\begin{eqnarray}\label{edbracketabstract}
[D_{i},Y(\psi,\x)]\phi=\partial_{x_{i}}Y(\psi,\x)\phi
\;\;\;\mbox{ for }\psi,\phi\in V,\; i=1,\dots, n.
\end{eqnarray}
\ep

\pf Using Corollary \ref{cassocdefprod} we obtain (\ref{edbracketabstract}) as
\begin{eqnarray}
& &([D_{i},Y(\psi,\z)]\phi)(\y)\nonumber\\
&=&D_{i}(Y(\psi,\z)\phi)(\y)-Y(\psi,\z)D_{i}\phi(\y)\nonumber\\
&=&\partial_{y_{i}}(Y(\psi,\z)\phi)(\y)-Y(\psi,\z)
D_{i}\phi(\y)\nonumber\\
&=&\partial_{y_{i}}\iota_{\0,\y,\z}(\psi(\z+\y)\phi (\y))
-\iota_{\0,\y,\z}(\psi(\z+\y)\partial_{y_{i}}\phi (\y))\nonumber\\
&=&\iota_{\0,\y,\z}\left(\left(\partial_{y_{i}}\psi(\z+\y)\right)
\phi (\y)\right)\nonumber\\
&=&\iota_{\0,\y,\z}\left(\left(\partial_{z_{i}}\psi(\z+\y)\right)
\phi (\y)\right)\nonumber\\
&=&\partial_{z_{i}}\iota_{\0,\y,\z}(\psi(\z+\y)\phi (\y))\nonumber\\
&=&\partial_{z_{i}}(Y(\psi,\z)\phi)(\y).\;\;\;\;\Box
\end{eqnarray}


The following result gives the locality of adjoint vertex operators
(vertex operators on the space of vertex operators):

\bp{padjointlocal}
Let $V$ be a vertex algebra of mutually local multi-variable 
vertex operators on $W$. Then
for $\psi,\phi\in V$, $Y(\psi,\x)$ and $Y(v,\y)$ 
are mutually local $G_{n}$-vertex operators on $V$.
\ep

\pf Given $\psi,\phi\in V$,
let $\theta\in V,\; w\in W$. Let $l\in \N$ be such that
$$\x^{l}\psi(\x)w,\;\; \x^{l}\phi(\x)w\in W[[\x]].$$
Then from Lemma \ref{ldefweakassociativity} for {\em any} $\rho\in V$,
\begin{eqnarray}
& &(\x+\z)^{l}(Y(\psi,\x)\rho)(\z)w
=(\x+\z)^{l}\psi(\x+\z)\rho(\z)w,\\
& &(\y+\z)^{l}(Y(\phi,\y)\theta)(\z)w
=(\y+\z)^{l}\phi(\y+\z)\theta(z)w.
\end{eqnarray}
Especially, this holds for $\rho=Y(\psi,\y)\theta\in V((\y))$.
Then
\begin{eqnarray}\label{eproductxy}
& &(\x+\z)^{l}(\y+\z)^{l}(Y(\psi,\x)Y(\phi,\y)\theta)(\z)w\nonumber\\
&=&(\x+\z)^{l}(\y+\z)^{l}\psi(\x+\z)(Y(\phi,\y)\theta)(\z)w\nonumber\\
&=&(\x+\z)^{l}(\y+\z)^{l}\psi(\x+\z)\phi(\y+\z)\theta(z)w.
\end{eqnarray}
Similarly, we have
\begin{eqnarray}\label{eproductyx}
(\x+\z)^{l}(\y+\z)^{l}(Y(\phi,\y)Y(\psi,\x)\theta)(\z)w
=(\x+\z)^{l}(\y+\z)^{l}\phi(\y+\z)\psi(\x+\z)\theta(z)w.
\end{eqnarray}
Let $k\in \N$ (only depending on $\psi$ and $\phi$) be such that 
\begin{eqnarray}\label{elocallocal}
(\x-\y)^{k}\psi(\x)\phi(\y)=(\x-\y)^{k}\phi(\y)\psi(\x).
\end{eqnarray}
Using (\ref{eproductxy}), (\ref{eproductyx}) and (\ref{elocallocal}) we get
\begin{eqnarray}
& &(\x-\y)^{k}(\x+\z)^{l}(\y+\z)^{l}(Y(\psi,\x)Y(\phi,\y)\theta)(\z)w
\nonumber\\
&=&(\x-\y)^{k}(\x+\z)^{l}(\y+\z)^{l}(Y(\phi,\y)Y(\psi,\x)\theta)(\z)w,
\end{eqnarray}
noting that $(\x-\y)^{k}=((\x+\z)-(\y+\z))^{k}$.
Because of the lower truncation condition, we may multiply 
both sides by $(\z+\x)^{-l}(\z+\y)^{-l}$ to get
\begin{eqnarray}
(\x-\y)^{k}(Y(\psi,\x)Y(\phi,\y)\theta)(\z)w
=(\x-\y)^{k}(Y(\phi,\y)Y(\psi,\x)\theta)(\z)w.
\end{eqnarray}
Note that $k$ does not depend on $w$. The assertion is 
proved. $\;\;\;\;\Box$

Combining Lemma \ref{ldoperatorvacuum},
Propositions \ref{pdoperatorproperty} and \ref{padjointlocal}
with Theorem \ref{tequivalentaxioms} and also using Proposition \ref{pjacobirep}
we immediately have (cf. [Li2]):

\bt{tmain}
Let $W$ be a vector space. Then any closed local subspace of $VO_{G_{n}}(W)$ 
is a vertex algebra with $W$ as a natural module.
In particular, any maximal local subspace of $VO_{G_{n}}(W)$ 
is a vertex algebra with $W$ as a natural module.
$\;\;\;\;\Box$
\et

Let $S$ be a set of pair-wise mutually local $G_{n}$-vertex operators 
on $W$. Recall from Remark \ref{rmaaximal} that
$\<S\>$ is  the smallest closed local subspace of $VO_{G_{n}}(W)$ that contains 
$S$ and the identity operator $1$ of $W$.
In view of Theorem \ref{tmain},  $\<S\>$ is a vertex algebra with $W$ 
as a natural module. Therefore we have (cf. [Li2]):

\bt{tmain2}
Let $W$ be a vector space. Then any set of pair-wise mutually local 
$G_{n}$-vertex operators on $W$ generates  
a canonical vertex algebra with $W$ as a natural module.$\;\;\;\;\Box$
\et

Combining Theorem \ref{tmain2} with Proposition \ref{paddderivation} 
we immediately have:

\bp{pgeneratingwithderivations}
Let $W$ be a vector space and let $\d=(d_{1},\dots,d_{n})$ be an $n$-tuple
of endomorphisms of $W$. Then any set of pair-wise mutually local 
$G_{n}$-vertex operators in $VO_{G_{n}}(W,\d)$ generates  
a canonical vertex algebra $V$ with $W$ as a natural module such that
\begin{eqnarray}
[d_{i},Y(v,\x)]=\partial_{x_{i}}Y(v,\x) \;\;\;\mbox{ for }i=1,\dots,n, \; v\in V.
\;\;\;\;\Box
\end{eqnarray}
\ep

The following is an analogue of a theorem of [FKRW] (see also [MP1-2], [X]):

\bp{pfkrw}
Let $V$ be a vector space, let $\d=(d_{1},\dots,d_{n})$ be an $n$-tuple
of endomorphisms of $V$, ${\bf 1}$ a vector in $V$ and
$U$ a subspace of $V$ equipped with a linear map
$Y$ from $U$ to $VO_{G_{n}}(V)$ such that\\
(C1) $V$ is generated from ${\bf 1}$ by all $u_{\m}$ for $u\in V,\;
\m\in \Z^{n}$, where $Y(u,\x)=\sum_{\m}u_{\m}\x^{-\m-1}$.\\
(C2) $Y(u,\x){\bf 1}\in V[[\x]]$ and $\lim_{\x\mapsto 0}Y(u,\x){\bf 1}=u$ 
for $u\in U$.\\
(C3) For $u\in U,\; 1\le i\le n$, $[d_{i},Y(u,\x)]=\partial_{x_{i}}Y(v,\x)$.\\
(C4) For $u,v\in U$, $Y(u,\x)$ and $Y(v,\y)$ are mutually local.\\
(C5) $d_{i}({\bf 1})=0$ for $i=1,\dots,n$.\\
Then there is a unique
$G_{n}$-vertex algebra structure on $V$ extending the linear map $Y$ on $U$
with ${\bf 1}$ being the vacuum vector.
\ep

\pf Due to (C1), the uniqueness follows immediately from the iterate
formula (\ref{eiterate}) and induction.
Set
\begin{eqnarray}
S=\{ Y(u,\x)\;|\; u\in U\}.
\end{eqnarray}
By (C3)-(C4), $S$ is a local set of $VO_{G_{n}}(V,\d)$. In view of 
Proposition \ref{pgeneratingwithderivations},
$S$ generates a $G_{n}$-vertex algebra $\tilde{V}$ with $V$ as a module
such that
\begin{eqnarray}\label{edbracketmoduleproof}
[d_{i},Y(\psi,\x)]=\partial_{x_{i}}Y(\psi,\x)\;\;\;\mbox{ for }i=1,\dots,
n,\; \psi\in \tilde{V}.
\end{eqnarray}
With (C2), in view of Lemma \ref{lvacuumlikeweak}, we have
\begin{eqnarray}\label{evacuumfkrw}
Y(\psi,\x){\bf 1}\in V[[\x]]\;\;\;\mbox{ for }\psi\in \tilde{V}.
\end{eqnarray}
and we have a 
$\tilde{V}$-homomorphism $f$ from $\tilde{V}$ to $V$ defined by
\begin{eqnarray}
f(\psi)=\lim_{\x\rightarrow 0}Y(\psi,\x){\bf 1}
=\lim_{\x\rightarrow 0}\psi(\x){\bf 1}
\;\;\;\mbox{ for }\psi\in \tilde{V}.
\end{eqnarray}
Note that we use $1$ (the identity map on $V$) for the vacuum vector of 
$\tilde{V}$ and ${\bf 1}$ is the given element of $V$.
Furthermore, by Proposition \ref{pisomorphism}, $f$ is an isomorphism.
Then $V$ has a vertex algebra structure 
transported from $\tilde{V}$ through the linear isomorphism $f$.
By (C2) we have $f(Y(u,\x))=u$ for $u\in U$.
Therefore, the vertex algebra structure on $V$ extends the linear 
map $Y$ on $U$. $\;\;\;\;\Box$

\end{document}